\documentclass[11pt]{pjm}
\usepackage{amssymb,amsmath}

\renewcommand{\labelenumi}{(\roman{enumi})}
\newcommand {\art}[6]{{\sc #1:} {#2.} {\em #3} {\bf #4} {(#5),} {#6.}}
\newcommand {\book}[5]{{\sc #1:} {``#2."} {#3,} {#4} {(#5).}}
\newcommand {\samp}[8]{{\sc #1:} {#2,} {{\em in:} ``#3,"} {(#4),}
                                  {#5,} {#6,} {(#7),} {#8.}}

\newcommand {\toappear}[3]{{\sc #1:} {#2.} {\em #3,} {\sl{to appear}.}}

\newcommand {\submitted}[2]{{\sc #1:} {#2,} {\sl submitted.}}

\newcommand \reals {[\omega ]^\omega }

\renewcommand \phi {\varphi}

\newcommand {\open}[1]{({#1})^{\omega}}
\newcommand {\nseg}[2]{({#1})^{#2*}}
\newcommand {\opf}[1]{({#1})^{\omega}_{\F}}

\newcommand {\nopeni}[3]{({#1})^{#2}_{#3}}

\newcommand {\res}{{\upharpoonright}}
\newcommand \parto {(\omega)^\omega}

\newcommand \comeager {({\rm m}^\st_1)}
\newcommand \meager {({\rm m}^\st_0)}

\newcommand \N {\operatorname{\mathbb N}}
\newcommand \NN {( \N)}

\newcommand \force {{\hspace{0.4mm}}{\rule{0.1mm}{2.4mm}}
                       {\hspace{0.5mm}}{\rule{0.1mm}{2.4mm}}
               {\rule [1.2mm]{2.3mm}{0.1mm}}{\hspace{0.4mm}}}

\newcommand \cp {\operatorname{cp}}
\newcommand \pc {\operatorname{pc}}
\newcommand \mdom {\operatorname{dom}}
\newcommand \mmin {\operatorname{min}}

\newcommand \mMin {\operatorname{Min}}
\newcommand \stem {\operatorname{stem}}
\newcommand \add {\operatorname{add}}
\newcommand \cov {\operatorname{cov}}

\newcommand \suc {\operatorname{succ}}

\newcommand \feq {\sqsupseteq}
\newcommand \ceq {\sqsubseteq}
\newcommand \aceq {\ceq^*}
\newcommand \aeq {\stackrel{*}{=}}
\newcommand \kap {\sqcap}

\newcommand \seg {\preccurlyeq}

\newcommand \F  {{\mathfrak F}}

\newcommand \filter {{\mathcal F}}
\newcommand \cP {{\mathcal P}}

\newcommand \cO {{\mathcal O}}
\newcommand \cJ {({\rm m}_0)}

\newcommand \cG {{\mathcal G}}
\newcommand \cu {{\operatorname{u}}}
\newcommand \fF {{\mathfrak F}}
\newcommand \fE {{\mathfrak E}}

\newcommand \fH {{\mathfrak H}}
\newcommand \fJ {\meager}
\newcommand \fh {{\mathfrak h}}

\newcommand \bG {\boldsymbol{G}}

\newcommand \BSigma {{\boldsymbol{\Sigma}}}
\newcommand \BPi {{\boldsymbol{\Pi}}}
\newcommand \BDelta {{\boldsymbol{\Delta}}}

\newcommand \st {{\scriptscriptstyle{\bigstar}}}
\newcommand \sst {{}^{\st}}
\newcommand \MM {{\mathbb M}}
\newcommand \MMs {{\mathbb M}^{\st}}
\newcommand {\MMU}[1] {{\mathbb M}_{#1}^{\st}}
\newcommand \Ms {M^{\st}}
\newcommand \Us {{\mathbb U}^{\st}}

\newcommand \M {{\mathbf M}}
\newcommand \bfN {{\mathbf N}}
\newcommand \V {{\mathbf V}}
\newcommand \bfL {{\mathbf L}}

\newcommand \La {{\mathbb L}}
\newcommand {\LU}[1]{{{\La}_{#1}^{\st}}}
\newcommand {\lU}[1]{{L_{#1}^{\st}}}
\newcommand {\Gen} {{\boldsymbol{X}_G}}
\newcommand {\bXG} {\Gen}
\newcommand {\gen} {X_G}

\newcommand \CH {{\sf CH}}
\newcommand \MA {{\sf MA}}
\newcommand \ZFC {{\sf ZFC}}
\newcommand \Ord {{\rm Ord}}

\newcommand \la {\langle}
\newcommand \ra {\rangle}
\newcommand \subs {\subseteq}

\theoremstyle{definition}
\newtheorem {nummer}{ }[section]
\newtheorem {thm}[nummer]{Theorem}
\newtheorem {conj}[nummer]{Conjecture}
\newtheorem {prop}[nummer]{Proposition}
\newtheorem {cor}[nummer]{Corollary}
\newtheorem {lm}[nummer]{Lemma}
\newtheorem {fct}[nummer]{Fact}
\newtheorem {obs}[nummer]{Observation}

\newtheorem {FPT}[nummer]{First Periodicity Theorem}
\newtheorem {SPT}[nummer]{Second Periodicity Theorem}
\newtheorem {plm}[nummer]{Preliminary Lemma}
\newtheorem {HJthm}[nummer]{Hales-Jewett Theorem}

\newcommand \poplm {\parindent=0pt
                   {\it{Proof of the Preliminary Lemma.}}\hspace*{1ex}
                    \parindent=2.2ex}

\def\qed{{\unskip\nobreak\hfil\penalty50\hskip8mm\hbox{}
  \nobreak\hfil
  {{\small{\bf{q.e.d.}}}}\parfillskip=0mm \par\smallskip}}

\begin{document}
\title[Approaching the
dual Ramsey property]{Techniques for approaching the dual
Ramsey property in the projective hierarchy}
\author[L.~Halbeisen, B.~L\"owe]{Lorenz Halbeisen, Benedikt L\"owe}
\thanks{The first author wishes to thank the {\sl
Swiss National Science Foundation\/} for supporting him. The second
author wishes to thank the {\sl Studienstiftung des deutschen Volkes}
for the funding of his travels to Boise ID and Las Vegas NV,
and in particular
Gunter Fuchs (Berlin).}
\subjclass{{\bf 03E15 03E40} 03E05 05A18 05D10 03E60}
\address{Department of Mathematics, University of California,
Evans Hall, Berkeley CA 94720, U.S.A.}
\email{halbeis@math.berkeley.edu}
\address{Mathematisches Institut, Rheinische
Friedrich--Wilhelms--Uni\-ver\-si\-t\"at Bonn, Beringstra{\ss}e 6,
53115 Bonn, Germany}
\email{loewe@math.uni-bonn.de}

\begin{abstract}
We define the dualizations of objects and concepts which are
essential for investigating the Ramsey property in the first levels
of the projective hierarchy, prove a forcing equivalence theorem for
dual Mathias forcing and dual Laver forcing, and show that the
Harrington-Kechris techniques for proving the Ramsey property from
determinacy work in the dualized case as well.
\end{abstract}
\maketitle

\section{INTRODUCTION}

Set Theory of the Reals is a subfield of Mathematical Logic mainly
concerned with the interplay between Forcing and Descriptive Set
Theory. One of the motivations behind Descriptive Set Theory is the
strong intuition that simple sets of real numbers should not display
irregular behaviour, or, in other words, they should be topologically
and measure theoretically nice.

In order to fill this statement with mathematical content, we should
make clear what we mean by ``simple'' and what we mean by ``nice''.
Both questions have a conventional and well known answer:
\begin{itemize}
\item
The measure of simplicity with which we categorize our sets of reals
is the projective hierarchy, in other words, the number of
quantifiers necessary to define the sets with a formula in first
order analysis (or second order arithmetic).
\item
A set should be considered ``nice" or ``regular" if it has the Baire
property in all naturally occurring topologies on the real numbers
and is a member of all natural conceivable $\sigma$-algebras.
\end{itemize}

Set Theory teaches us that the axioms of {\sf ZFC} do not entail a
formal version of these intuitions: It is consistent with {\sf ZFC}
that there are irregular sets already at the first level
of the projective hierarchy.\footnote{In G\"odel's
Constructible Universe $\bfL$ there is a $\BDelta^1_2$ set which is not
Lebesgue measurable and which does not have the Baire property. Worse
still, there is an uncountable $\boldsymbol{\Pi}^1_1$ set with no
perfect subset and a $\BPi^1_1$ set which is not Martin measurable.}
Thus the focus shifts from proving that all simple sets are nice to
investigating the situations under which our intuitions are met by
the facts.

A whole array of research in this direction is dealing with the
second level of the projective hierarchy. Solovay provided us
with the prototype of a characterization theorem for the second
level:

\begin{thm}
The following are equivalent:
\begin{enumerate}
\item Every $\BSigma^1_2$ set of reals has the Baire property.
\item For every real $a\in\reals$ the set of Cohen generic reals over the
model ${\bfL}[a]$ is comeager in the standard topology on the real
numbers.
\end{enumerate}
\end{thm}

One could call a theorem like this a ``transcendence principle over
the constructible universe''. These principles connect the theory of
forcing and the topological properties of the reals. Comparable
theorems have been proved in \cite{JudahShelah.12} (for the
$\BDelta^1_2$ level) and in \cite{BrLoe99} (for different topologies
and $\sigma$-algebras).

A particularly interesting instance of niceness in the above sense is
the {\bf Ramsey property}, a topological property which is deeply
connected to Ramsey theory and infinitary combinatorics. The Ramsey
property is linked to a forcing notion called {\bf Mathias forcing},
introduced by Mathias in \cite{Mathias}, and Judah and Shelah were
able to obtain the following Solovay--type characterization for it
({\sl cf.}~\cite[Theorem\,2.7 \& Theorem\,2.8]{JudahShelah.12}):

\begin{thm}\label{JuShRamsey}
The following are equivalent:
\begin{enumerate}
\item Every $\boldsymbol{\Sigma}^1_2$ set of reals has the Ramsey property.
\item Every $\boldsymbol{\Delta}^1_2$ set of reals has the Ramsey property.
\item For every real $a\in\reals$ the set $\{r\in\reals : r$ is
Ramsey over $\bfL[a][\filter^r]\}$\footnote{A real $r$ is Ramsey over
$\bfL [a][\filter^r]$ if and only if $\filter^r:=D_r\cap\bfL
[a][D_r]$ forms an ultrafilter in $\bfL [a][\filter^r]$, where
$D_r:=\{r'\in\reals
:r\subs^* r'\}$, and $r$ is $\mathbb{L}_{\filter^r}$-generic over $\bfL
[a][\filter^r]$, where $\mathbb{L}_{\filter^r}$ is Laver forcing
restricted to $\filter^r$.} is comeager in the Ellentuck topology.
\end{enumerate}
\end{thm}

One connection to Mathias forcing is given by the following result
({\sl cf.}~\cite[Theorem\,4.1]{LoriJudah}):

\begin{prop} If $\bfN$ is any model of \ZFC, then the following are equivalent:
\begin{enumerate}
\item $\bfN$ is a model in which every $\BSigma^1_2$ set is
Ramsey, and
\item $\bfN$ is
$\BSigma^1_3$-Mathias-absolute.\footnote{Similar characterizations
exist also for some other properties, {\sl e.g.} for Lebesgue
measurability and Baire property ({\sl
cf.}~\cite[Theorem\,9.3.8]{BartoszynskiJudah}).}
\end{enumerate}
\end{prop}

As the Ramsey property talks about infinite subsets of the natural
numbers, it is easily dualized by something we shall call the {\bf
dual Ramsey property}, talking about infinite partitions of the
natural numbers.\footnote{Infinite subsets can be seen as images of
injective functions and infinite partitions can be seen as preimages
of surjective functions, so the move from infinite subsets to
infinite partitions actually is a dualization process.} This property
has been introduced by Carlson and Simpson in \cite{CarlsonSimpson}
and further investigated in \cite{Lorisha} and \cite{Lorisym}.

One thing that is striking about the relationship between the
Ramsey property and the dual
Ramsey property are the distinctive symmetries and asymmetries. This
paper can be understood as a catalogue of some of the similarities;
in fact, one could see parts of this paper as an attempt to reach the
obvious dualization of Theorem \ref{JuShRamsey}:

\begin{conj}\label{OurConjecture}
The following are equivalent:
\begin{enumerate}
\item Every $\BSigma^1_2$ set of reals has the dual Ramsey property.
\item Every $\BDelta^1_2$ set of reals has the dual Ramsey property.
\item For every real $a\in\reals$ the set $\{R : R$ is dual Ramsey over
$\bfL[a][D^R]\}$ is comeager in the dual Ellentuck topology.
\end{enumerate}
\end{conj}

In order to approach this conjecture and to give an idea what ``dual
Ramsey over $\bfL[a][D^R]$'' could mean, several of the techniques of
\cite{JudahShelah.12} and \cite{Mathias} have to be adapted to the
new environment:

Mathias forcing has a characteristic product form $\MM =
{\cP}(\omega)/{\rm fin} * {\MM}_{\mathbf U}$ where $\mathbf U$ is the
canonical name for the generic ultrafilter added by
${\cP}(\omega)/{\rm fin}$. This ultrafilter is in fact a Ramsey
ultrafilter,\footnote{\label{fn:ramseyuf}A
set $\filter\subs\reals$ is a Ramsey filter
if $\filter$ is a filter and for any colouring $\tau:[\omega]^n\to
r+1$ (with $n,r\in\omega$) there is an $x\in\filter$ such that
$\tau\res [x]^n$ is constant. Notice that every Ramsey filter is an
ultrafilter.} and Judah and Shelah show in their
\cite{JudahShelah.12} that Mathias forcing relative to an ultrafilter
is forcing-equivalent to Laver forcing relative to the same
ultrafilter, provided that the ultrafilter is Ramsey ({\sl
cf.}~\cite[Theorem\,1.20\,(i)]{JudahShelah.12}):
\begin{thm}\label{LavMath}
Let $\filter$ be a Ramsey ultrafilter. Then the forcing notions
${\mathbb L}_{\filter}$ and ${\MM}_{\filter}$ are equivalent.
\end{thm}

This theorem was our motivation to search for a dual version of Laver
forcing and the dualization of Ramsey ultrafilters to work towards a
dualization of Theorem~\ref{JuShRamsey}.

\bigskip

In our dualized situation there are many things to be done to make
sense of the dualized versions: One has to find a dualized version of
${\cP}(\omega)/{\rm fin}$ and to prove the corresponding product form
of dual Mathias forcing (already done in \cite{Lorisym}), one has to
find a dualized version of Ramsey ultrafilters, and one has to make
explicit what Laver forcing in this context is supposed to mean.

Section~\ref{sec:defnot} of this paper defines all the dualized
notions needed for the technical work on the dual Ramsey property. In
Section~\ref{sec:druf}, the reader will find a couple of facts about
a dualization of Ramsey ultrafilters; their connection to the game
filters from \cite{Lorisym} is given in the appendix.
Section~\ref{sec:Laver*} moves on to discuss dual Laver forcing and
proves the dualized version of Theorem \ref{LavMath}.

In Sections \ref{sec:simple} and \ref{sec:det} we investigate the extent
of sets with the dual Ramsey property in the projective hierarchy. In
Section~\ref{sec:simple} we prove a couple of consistency results
for the first three levels of the projective hierarchy.
After that,
Section~\ref{sec:det} looks at the dual Ramsey property from a
completely different angle: If we assume an appropriate amount of
determinacy, we know that a large collection of sets has the Ramsey
property. This result is not at all immediate from the Banach-Mazur
game for the topology associated with the Ramsey
property.\footnote{The obstacle is that playing basic open sets in
this topology cannot be coded by natural numbers. So the Banach-Mazur
games essentially needs determinacy for games with real moves, {\sl e.g.}, ${\sf
PD}_{\mathbb R}$. This
is connected to the famous open question whether {\sf AD} implies
that every set has the Ramsey property ({\sl cf.} \cite[Question
27.18]{Kan94}).} However, in that section we note that the
Harrington-Kechris technique of proving the Ramsey property from
standard determinacy ({\sl cf.}~\cite{HarringtonKechris}) alone works
for the dualized case as well.

\medskip

It should be mentioned that the technicalities of the dualization
process are not always as easy as they seem in retrospect. Finding
the correct and natural dualizations for the interesting notions from
the classical case is the most challenging part in this project.
After the right dualizations are at hand, in most cases one can
follow the classical proofs. So, the merits of this paper lie mainly
in the definitions that make the proofs nice and easy and give a
proper and firmly rooted understanding of the symmetries. This is
also the reason for the unproportional size of Section
\ref{sec:defnot} compared to the other sections.

\section{DEFINITIONS AND NOTATIONS}\label{sec:defnot}
{\setcounter{subsection}{-1}}

\subsection{Set-theoretic notation}

Most of our set-theoretic notation is standard and can be found in
textbooks like \cite{Jechbook}, \cite{Kunen} or
\cite{BartoszynskiJudah}. For the definitions and some basic facts
concerning the projective hierarchy we refer the reader to
\cite[\S 12]{Kan94}.

We shall consider the set $\reals$ as the set of real numbers. For the
Turing join of two reals $x$ and $y$ ({\sl i.e.}, coding two reals
into one), we use the standard notation $x\oplus y$.

\subsection{Partitions}

A set $P\subs{\mathcal P}(S)$ is a partition of the set $S$ if
$\emptyset\notin A$, $\bigcup P=S$ and for all distinct
$p_1,p_2\in P$ we have $p_1\cap p_2=\emptyset$. An element of a
partition $P$ is also called a {\bf block\/} of $P$ and $\mdom
(P):=\bigcup P$ is called the {\bf domain\/} of $P$. A partition
$P$ is called infinite, if $|P|$ is infinite, where $|P|$ denotes
the cardinality of the set $P$. The equivalence relation on $S$
uniquely determined by
a partition $P$ is denoted by $\sim_P$.

\smallskip

Let $P$ and $Q$ be two arbitrary partitions. We say that $P$ is {\bf
coarser\/} than $Q$ (or that $Q$ is {\bf finer\/} than $P$) and write
$P\ceq Q$, if for all blocks $p\in P$, the set $p\cap\mdom (Q)$ is
the union of some sets $q_i\cap\mdom (P)$, where each $q_i$ is a
block of $Q$. Let $P\kap Q$ be the finest partition which is coarser
than $P$ and $Q$ with $\mdom (P\kap Q)= \mdom (P)\cup\mdom (Q)$. We
say that $P$ is {\bf almost coarser\/} than $Q$ and write $P\aceq Q$
if there is a partition $R$ such that $\mdom (R)$ is finite and
$R\kap P\ceq Q$. If $P\aceq Q$ and $Q\aceq P$, then we write $P\aeq
Q$.\footnote{We choose this notation because the properties of $\ceq$
and $\kap$ are similar to those of $\subs$ and $\cap$.}


Let $P$ and $Q$ be two partitions. If for each $p\in P$ there is
a $q\in Q$ such that $p=q\cap\mdom (P)$, we write $P\seg Q$. Note
that $P\seg Q$ implies $\mdom (P)\subs\mdom (Q)$.

For $x\subs\omega$ let $\mmin (x):=\bigcap x$. If $P$ is a partition
with $\mdom (P)\subs\omega$, then $\mMin (P):=\{\mmin (p):p\in P\}$;
and for $n\in\omega$, $P(n)$ denotes the unique block $p\in P$ such
that $|\mmin (p)\cap\mMin (P)|=n+1$.

\medskip

The set of all infinite partitions of $\omega$ is denoted by
$\parto$; and the set of all partitions $s$ with $\mdom (s)\in\omega$
is denoted by $\NN$.

For $s\in\NN$, let $s^*$ denote the partition $s\cup\big{\{}\{\mdom
(s)\}\big{\}}$. Notice that $|s^*|=|s|+1$.

For a natural number $n$, let $\nseg{\omega}{n}$ denote the set of
all $u\in\NN$ such that $|u|=n$. Further, for $n\in\omega$ and
$X\in\parto$ let $$\nseg{X}{n}:=\{u\in\NN
:|u|=n\wedge u^*\ceq X\}\,,$$ and for $s\in\NN$ such that $|s|\le n$
and $s\ceq X$, let $$\nseg{s,X}{n}:=\{u\in\NN
:|u|=n\wedge s\seg u\wedge u^*\ceq X\}\,.$$

It will be convenient to consider $\omega$ as the partition which
contains only singletons, and therefore, for $s\in\NN$,
$\nseg{s,\omega}{n}:=\{u\in\NN :|u|=n\wedge s\seg u\}$.

\bigskip

A family $\F\subs\parto$ is called a {\bf filter} if
\begin{itemize}
\item[{($\alpha$)}] $\emptyset\notin\F$;
\item[{($\beta$)}] If $X\in\F$ and $X\ceq Y$, then $Y\in\F$;
\item[{($\gamma$)}] If $X$ and $Y$ belong to $\F$, then $X\kap Y\in\F$.
\end{itemize}
Further, we call $\F\subs\parto$ an {\bf ultrafilter} if $\F$ is a
filter which is not properly contained in any filter. Notice that if
$X\in\F$ and $\F\subs\parto$ is an ultrafilter, then each
$Y\in\parto$ with $Y\aeq X$ belongs to $\F$, too.

\subsection{The dual Ellentuck topology and the dual Ramsey property}

Let $X\in\parto$ and $s\in\NN$ be such that $s\ceq X$. Then
$${\open{s,X}}:=\{Y\in\parto :s\seg Y\ceq X\}$$ and $$\open{X}
:={\open{\emptyset,X}}=\{Y\in\parto: Y\ceq X\}\,.$$

Obviously, this definition depends on the model we are working in,
so, if this should become important, we denote by $\open{s,X}_\bfN$
the corres\-ponding set interpreted in the model $\bfN$.

\smallskip

Let the basic open sets on $\parto$ be $\emptyset$ and the sets
${\open{s,X}}$, where $s$ and $X$ are as above. These sets are called
the {\bf dual Ellentuck neighbourhoods}. The topology induced by the
dual Ellentuck neighbourhoods is called the {\bf dual Ellentuck
topology\/} ({\sl cf.}~\cite{CarlsonSimpson}).

\medskip

A family $A\subs\parto$ has the {\bf dual Ramsey property\/} (or just
{\bf is dual Ramsey}) if and only if there is a partition $X\in\parto$ such that
either $\open{X}\subseteq A$ or $\open{X}\cap A = \emptyset$.

\medskip

Closely related (but stronger) is the notion of a completely dual Ramsey set: A set
$A\subs\parto$ is said to be {\bf completely dual Ramsey} if and only
if for each dual Ellentuck neighbourhood $\open{s,X}$ there is a
$Y\in\open{s,X}$ such that $\open{s,Y}\subs A$ or $\open{s,Y}\cap
A=\emptyset$. If we are always in the latter case, then $A$ is called
{\bf completely dual Ramsey-null}. Concerning projective sets it is
not clear if ``completely dual Ramsey'' is really stronger than just
``dual Ramsey'', because we cannot simply translate Lemma~2.1 of
\cite{BrLoe99}, where it is shown under other things that ``Ramsey''
and ``completely Ramsey'' coincide with respect to projective sets.

\smallskip

Carlson and Simpson proved in~\cite{CarlsonSimpson}
that a set $A$ is completely dual Ramsey if and only if $A$ has the Baire
property with respect to the dual Ellentuck topology and $A$ is completely dual
Ramsey-null if and only if $A$ is meager with respect to the dual
Ellentuck topology.\footnote{A set $S$ has the Baire property if there is
an Borel set $B$ such that the symmetric difference $S\triangle B$ is
meager, where a meager set is the union of countably many nowhere
dense sets.} As a matter of fact we like to mention that in the dual
Ellentuck topology every meager set is nowhere dense and hence, the
dual Ellentuck topology is a Baire topology ({\sl i.e.}, no open set
is meager). This corresponds to the similar facts about ``being
completely Ramsey'' and the Ellentuck topology.\footnote{{\sl Cf.}
\cite{Ellentuck} \& \cite[\S 19.D]{Kechris}.}

\subsection{Dual Mathias forcing}\label{subsec:DMF}

The conditions of the {\bf dual Mathias forcing\/} $\MMs=\la
\Ms,\le\ra$ are the pairs $\la s,X\ra$ such that $\open{s,X}$ is a
non-empty dual Ellentuck neighbourhood, and the partial order is defined by

$$\langle s,X\rangle\leq\langle t,Y\rangle\ \Leftrightarrow\
{\open{s,X}}\subs{\open{t,Y}}\,.$$

If $\langle s,X\rangle$ is an
$\MMs$-condition, then we call $s$ the {\bf{stem}\/} of the
condition.

\medskip

If $G$ is $\MMs$-generic over $\bfN$, then $G$ induces in a canonical
way an infinite partition $\gen\in\parto$ such that
${\bfN}[G]={\bfN}[\gen]$, and therefore we consider the partition
$\gen$ as the generic object. Thus, we can reconstruct the original
$G$ from $\gen$ by observing that

$$\la s,X\ra\in G ~\iff~ \gen\in\open{s,X}_{\bfN[G]}\,.$$

\smallskip

Since the dual Ellentuck topology is innately connected with dual
Mathias forcing, we choose the following notation for meager and
comeager sets in the dual Ellentuck topology:

\begin{eqnarray*}
A\in\meager & \iff & A\mbox{ is dual Ellentuck meager, and}\\
A\in\comeager & \iff & A\mbox{ is dual Ellentuck comeager,}\\
&&\mbox{{\sl i.e.}, } \parto\setminus A\mbox{ is dual Ellentuck meager.}
\end{eqnarray*}

Since $A$ is dual Ellentuck meager if and only if $A$ is completely
dual Ramsey-null, $\meager\subs{\mathcal{P}}(\parto)$ is also the
ideal of completely dual Ramsey-null sets.

\medskip

The following fact gives two properties of dual Mathias forcing which
also hold for Mathias forcing.

\begin{fct} \label{HomPureDec}
If $\gen$ is $\MMs$-generic and $Y\in\open{\gen}$,
then $Y$ is $\MMs$-generic as well (we will call this property the
{\bf homogeneity property}); and therefore, dual Mathias forcing
is proper. Moreover, for any sentence $\Phi$ of the forcing language
$\MMs$ and for any $\MMs$-condition $\la s,X\ra$, there is an
$\MMs$-condition $\la s,Y\ra\le\la s,X\ra$ such that $\la
s,Y\ra\force_{\MMs}\Phi$ or $\la s,Y\ra\force_{\MMs}\neg\Phi$ (this
property is called {\bf pure decision}).
\end{fct}

\proof For a proof, {\sl cf.}
\cite[Theorem 5.5 \& Theorem 5.2]{CarlsonSimpson}.\qed

\smallskip

As an immediate
consequence we get that the set of dual Mathias generic partitions
over every model $\bfN$ is either empty or a non--meager set
which is completely dual Ramsey.

\medskip

Like Mathias forcing, dual Mathias forcing has also a
characteristic product form.

\smallskip

Let $\Us=\la\parto,\le\ra$ be the partial order defined as follows:
$$X\leq Y \ \Leftrightarrow\ X\ceq^* Y\,.$$

$\Us$ is the natural dualization of $\cP(\omega)/{\rm fin}$.

\smallskip

For a family $\fE\subs\parto$ we define the {\bf restricted dual
Mathias forcing} $\MMU{\fE}$ as follows. The conditions of
$\MMU{\fE}=\la M^{\st}_{\fE},\le\ra$ are the $\MMs$-conditions $\la
s,X\ra$ such that $X\in\fE$.

Now we get

\begin{fct} $\MMs=\Us * \MMs_{\bG}$, where
${\bG}$ is the canonical name for the $\Us$-generic object.
\end{fct}

\proof For a proof, {\sl cf.} \cite[Fact\,2.5]{Lorisym}. \qed

\subsection{Restricted dual Laver forcing}

In order to define the forcing notion which will be investigated
later on, we first have to give some notations.

For $T\subs\NN$ and $t\in T$ we define the {\bf successor set of $t$
in $T$} as follows: $$\suc_T(t):=\{u\in T:t\seg u\wedge
|u|=|t|+1\}\,.$$

Let $\fE\subs\parto$ be any non-empty family (later on we investigate
only the case when $\fE$ is an ultrafilter).

With respect to $\fE$, we define the {\bf dual Laver forcing
restricted to $\fE$}, denoted by $\LU{\fE}=\la \lU{\fE},\le\ra$, as
follows:

\begin{itemize}
\item[{($\alpha$)}] $p\in\lU{\fE}$ if and only if $p\subs\NN$ with the
property that there is an $s\in p$ (denoted $\stem(p)$) such that for
all $t\in p$ we have $s\seg t$.

\item[{($\beta$)}] There exists a set $\{X_t^p:t\in p\}\subs\fE$ such that
for $t\in p$ we have $t^*\ceq X_t^p$ and
$$\suc_p(t)=\{u:u\in\nseg{t^*,X_t^p}{(|t|+1)}\}\,.$$ Further, for
$t,u\in p$ with $t\seg u$ we have
$$\open{u,X_u^p}\subs\open{t,X_t^p}\,,$$ and if $\mdom(t)=\mdom(u)$
and $t\ceq u$, then $$X_t^p=X_u^p\,.$$

\item[{($\gamma$)}] For two $\LU{\fE}$-conditions $p$ and $q$ we stipulate
$$p\le q\ \iff\ p\subs q\,.$$
\end{itemize}

Notice that $p\le q$ implies $\stem(q)\seg\stem(p)$ and hence, if
$G\subs \lU{\fE}$ is $\LU{\fE}$-generic over some ${\bfN}$, then the
set $\{s:s=\stem(p)\;\text{for some $p\in G$}\}$ forms in a canonical
way a partition $\gen\in\parto$. Moreover, ${\bfN}[G]={\bfN}[\gen]$
and therefore we may consider also the partition $\gen$ as the
$\LU{\fE}$-generic object.

For an $\LU{\fE}$-condition $p$ we call a partition $X\in\parto$ a
{\bf branch of} $p$ if each $t\in\NN$ with $t^*\ceq X$ belongs to
$p$.

\begin{fct}\label{fct:sub-branch}
If $X$ is a branch of the $\LU{\fE}$-condition $p$ where $\stem(p)=s$
and $Y\in\open{s,X}$, then $Y$ is a branch of $p$, too.
\end{fct}

\proof This follows immediately from ($\beta$). \qed

\subsection{Special ultrafilters on $\boldsymbol{\parto}$}

A family $\F$ has the {\bf segment colouring property} (or just {\bf
scp}) if for any $s\ceq X\in\F$ with $|s|=n$ and for any colouring
$\pi :\nseg{s,X}{(n+k)}\to r$, where $r$ and $n+k$ are positive
natural numbers, there is a $Y\in\open{s,X}\cap\F$ such that
$\nseg{s,Y}{(n+k)}$ is monochromatic.

A family $\F\subs\parto$ is an {\bf scp-filter} if $\F$ is a filter
which has the segment colouring property.

\smallskip

In Section \ref{sec:appendix} we shall introduce the notion of game
filters (from \cite{Lorisym}) and show that game filters are
scp-filters.

\begin{fct}
If $\F\subs\parto$ is an scp-filter, then $\F$ is an ultrafilter.
\end{fct}

\proof Let $\F\subs\parto$ be an scp-filter and assume that there
exists an $X\in\parto$ such that for every $Y\in\F$, $X\kap
Y\in\parto$. Let $\pi:\nseg{\omega}{n}\to 2$ be such that $\pi (s)=0$
if $s\in\nseg{X}{n}$, otherwise $\pi (s)=1$. Because $\F$ has the
segment colouring property, we find a $Y\in\F$ such that
$\pi\res{\nseg{Y}{n}}$ is constant. If $\pi\res{\nseg{Y}{n}}=\{1\}$,
then $X\kap Y\notin\parto$ which contradicts the assumption. Thus,
$\pi\res{\nseg{Y}{n}}=\{0\}$, which implies $X\in\F$ and hence, the
filter $\F$ is maximal. \qed

\medskip

A family $\F$ is {\bf diagonalizable} if for any $\LU{\F}$-condition
$p$, there is a partition $X\in\F$ such that $X$ is a branch of $p$.
Notice that a diagonalizable family can also be characterized by a
two player game, where the $\LU{\F}$-condition $p$ can be considered
as a strategy for player~I.

\medskip

A family $\F$ is a {\bf Ramsey$\sst$ filter} if $\F$ is a
diagonalizable scp-filter.

\medskip

In Footnote \ref{fn:ramseyuf}
we have defined Ramsey ultrafilters over $\omega$ in terms of
colourings. This definition corresponds to the definition of
scp-filters. On the other hand, Galvin and Shelah proved that Ramsey
ultrafilters can be characterized as well by a two player game
without a winning strategy for player~I, where a winning strategy for
player~I is in fact a restricted Laver-condition ({\sl
cf.}~\cite[Theorem\,4.5.3]{BartoszynskiJudah}). This definition of
Ramsey ultrafilters corresponds to diagonalizable filters.
It is possible that the notions of ``scp-filters'' and
``diagonalizable filters'' are equivalent, but this is still open.

\smallskip

Beyond the dualization of the notion of a Ramsey ultrafilter, the
dualization process leading from $\reals$ to $\parto$ has interesting
consequences for the spaces of ultrafilters on these spaces. These
consequences belong to the asymmetrical aspects of the relationship
between $\reals$ and $\parto$ and are the point of focus in
\cite{luxor}.

\subsection{Switching between reals and partitions}

We fix $\flat : [\omega]^2\to\omega$ to be any arithmetic bijection
between the set of pairs of natural numbers and $\omega$.

Let $x\in \reals$; then the set $\mbox{trans}(x)\subseteq\omega$ is
defined by
\begin{eqnarray*}
n\in{\rm trans}(x)&:\iff& \exists
s\in\omega^{<\omega}~\Big(~n=\flat(s(0),s(|s|-1))\mbox{ and }\\
&&\forall k\in |s|-1 (\flat(s(k),s(k+1))\in x)~\Big).
\end{eqnarray*}

As the name suggests, ${\rm trans}(x)$ is the set of codes of pairs
in the transitive closure of the relation $\flat(k,\ell)\in x$. A
real $x$ is called {\bf transitive} if ${\rm trans}(x) = x$.

Note that in general $\mbox{trans}(x)\subseteq x$ and that the
relation $$R_x(k,\ell):\iff\flat(k,\ell)\in{\rm trans}(x)$$ is
symmetric (by choice of the domain of $\flat$) and transitive. Thus,
if $x\in [\omega]^\omega$, we can consider $x$ as a partition (by
reflexivization of $R_x$) via

$$n\sim_x m~:\iff~n=m\mbox{ or }\flat(n,m)\in{\rm trans}(x).$$

We call this partition the {\bf corresponding partition} of $x\in
\reals$, and denote it by $\cp (x)$. Note that $\cp (x)\in\parto$ if

$$\forall k\exists n>k\forall m<n(\neg(n \sim_x m))$$

and further if $y\subseteq x$, then $\cp (y)\feq\cp (x)$.

We encode a partition $X$ of $\omega$ by a real $\pc (X)$ (the {\bf
partition code\/} of $X$) as follows. $$\pc (X):=\{k\in\omega
:\exists n \exists m(k=\flat(n,m)\wedge (n\sim_{X} m))\}.$$ Note
that if $X\ceq Y$ then $\pc (X)\supseteq\pc (Y)$.

Notice that both the function $\pc$ and the function $\cp$ are
arithmetic, and that they are in a sense inverse to each other:

\begin{obs}\label{INVERSE}
For every $X\in\parto$ and every $x\in\reals$ the following hold:
\begin{enumerate}
\item $\cp(\pc(X)) = X$, and
\item if $x$ is transitive, then $\pc(\cp(x)) = x$.
\end{enumerate}
\end{obs}

\medskip

Now, a set $A\subs\reals$ has the {\bf dual Ramsey property\/} (or
just {\bf is dual Ramsey\/}) if and only if the set $\{X\in\parto:
\exists x\in A(X=\cp (x))\}$ has the dual Ramsey property. By
Observation \ref{INVERSE}, this is equivalent to saying that the set
$\{ X\in\parto : \pc(X)\in A\}$ has the dual Ramsey property.

\smallskip

By the definition of the dual Ramsey property we have that every
$\BSigma^1_n$ set is dual Ramsey if and only if every $\BPi^1_n$ set
is dual Ramsey. Further we have by \cite[Lemma 7.2]{Lorisym} that if
every $\BSigma^1_n$ set is dual Ramsey then every $\BSigma^1_n$ set
has the classical Ramsey property.

\medskip

As a matter of fact we like to mention the following

\begin{prop}\label{prop:Delta.1n}
If every $\BDelta^1_n$ set has the dual Ramsey property, then every
$\BDelta^1_n$ set has the Ramsey property.
\end{prop}

\proof Suppose $A$ is a $\BDelta^1_n$ set of reals. Let $\phi$ be a
$\BSigma^1_n$ formula and $\psi$ be a $\BPi^1_n$ formula witnessing
this, {\sl i.e.}, $$x \in A\ \;\iff\ \;\phi(x)\ \;\iff\
\;\psi(x).$$ To show that $A$ is Ramsey we define a different
$\BDelta^1_n$ set by formulae $\phi^*$ and $\psi^*$ as follows:
\begin{align*} \phi^*(v) & :\iff \exists w ( w = \mMin(\cp(v))
\wedge \phi(w))\\ \psi^*(v) & :\iff \forall w ( w = \mMin(\cp(v))
\rightarrow \psi(w))
\end{align*}

Obviously, $\phi^*$ is $\BSigma^1_n$ and $\psi^*$ is $\BPi^1_n$, and
since $\mMin(\cp(v))$ is uniquely determined for each $v$, these two
formulae are equivalent and hence define a $\BDelta^1_n$ set $A^*$ of
reals. The rest of the proof is exactly as in
\cite[Lemma\,7.2]{Lorisym}. \qed

And as a corollary we get

\begin{cor}
If every $\BDelta^1_2$ set is dual Ramsey, then every $\BSigma^1_2$
set is Ramsey.
\end{cor}

\proof This follows immediately from Proposition~\ref{prop:Delta.1n}
by Theorem~\ref{JuShRamsey}. \qed

\subsection{Descriptive Set Theory of the Cabal}
\label{subsec:cabal}

For our results in Section~\ref{sec:det} we shall need some basic
notions of the Descriptive Set Theory of the Cabal Seminar.
Everything we lay out here can be found in \cite{Moschovakis}, our
account is just for the convenience of the more combinatorially
oriented reader who might be unfamiliar with the language of the
Cabal.

We shall presuppose basic knowledge with the standard notation for
determinacy and the elementary results of the theory of perfect
information games as outlined in \cite[\S 27]{Kan94}.

\medskip

Let $X$ be a set of reals and $\alpha\in\Ord$. Any surjective
function $\phi: X\to\alpha$ is called a {\bf norm on} $X$. The
ordinal $\alpha$ is called the {\bf length of $\phi$}. A family $\Phi
:= \la \phi_n : n\in\omega\ra$ of norms on $X$ is called a {\bf scale
on} $X$ if for every sequence $\la x_i : i\in\omega\ra\subseteq X$
and every $n\in\omega$ the following holds: If $\la\phi_n(x_i) :
i\in\omega\ra$ is eventually constant, say, equal to $\lambda_n$,
then $x:= \lim_{i\in\omega} x_i \in X$ and
$\phi_n(x)\leq\lambda_n$.\footnote{For the basic theory of scales,
{\sl cf.} \cite{KechrisMoschovakis}.}

Let $\Gamma$ be any pointclass, $\phi$ any norm on $X$, and $\Phi$
any scale on $X$. We shall call $\phi$ a $\Gamma$ {\bf norm} if there
are two relations $R$ and $R^*$ in $\Gamma$ such that:

$$y\in X \Rightarrow  \forall x \big( (x\in X \wedge
\phi(x)\leq\phi(y)) \iff R(x,y) \iff \neg R^*(x,y)\big).$$

\medskip

We call a scale $\Phi$ a $\Gamma$ {\bf scale} if all norms $\phi_n$
occurring in $\Phi$ are $\Gamma$ norms, uniformly in $n$.\footnote{A
more precise definition can be found in \cite[p.\,228]{Moschovakis}.}
We shall say that a set $X$ {\bf admits a $\Gamma$ norm} ({\bf a
$\Gamma$ scale}) if there is a norm (a scale) on $X$ that is a
$\Gamma$ norm (a $\Gamma$ scale).

The fundamental theorems connecting determinacy, norms and scales are
the ``Periodicity Theorems'' of \cite{AdMo68}, \cite{Ma68} and
\cite{Mo71}. In the following we shall need the first two Periodicity
Theorems in special cases:

\begin{FPT}\label{FPT}
Suppose that ${\rm Det}(\boldsymbol{\Delta}^1_{2n})$ holds and
$x\in\reals$ is a real. Then every $\Pi^1_{2n+1}(x)$ set admits a
$\Pi^1_{2n+1}(x)$ norm.
\end{FPT}

\begin{SPT}\label{SPT}
Suppose that ${\rm Det}(\boldsymbol{\Delta}^1_{2n})$ holds and
$x\in\reals$ is a real. Then every $\Pi^1_{2n+1}(x)$ set admits a
$\Pi^1_{2n+1}(x)$ scale.
\end{SPT}

For proofs of these theorems (in a much more general formulation), we
refer the reader to \cite[6B.1 \& 6C.3]{Moschovakis}.

We define (for every $n\in\omega$) the projective ordinals by
$$\boldsymbol{\delta}^1_{n} := \sup\{\|{\leq}\| : {\leq} \mbox{ is a
$\Delta^1_n$ prewellordering on }\reals\},$$

and note that for every $\Pi^1_{2n+1}$ complete set the length of every
$\Pi^1_{2n+1}$ norm on it is
exactly $\boldsymbol{\delta}^1_{2n+1}$ (\cite[4C.14]{Moschovakis}).

\bigskip

Let
$P^x_{2n+1}$ be a $\Pi^1_{2n+1}(x)$ complete set of reals. Assuming
${\rm Det}(\BDelta^1_{2n})$ we get a $\Pi^1_{2n+1}(x)$ scale
$\Phi^x = \{\varphi^x_m: m\in\omega\}$ for $P^x_{2n+1}$ by Theorem \ref{SPT}.

For any real $y\in P^x_{2n+1}$, we denote by $\Phi^x(y)$ the sequence
of ordinals determined by the scale, {\sl i.e.}, $\Phi^x(y) = \langle
\varphi^x_m(y) : m\in\omega\rangle$.

\smallskip

Now let
$$T^x_{2n+1} := \{ \langle y\res m, \Phi^x(y)\res m\rangle : y\in
P^x_{2n+1}, m\in\omega\}$$

be the {\bf tree associated to} $\Phi^x$. By the remark about the
lengths of norms, it is a tree on
$\omega\times\boldsymbol{\delta}^1_{2n+1}$. If $x$ is any recursive
real, we write $T_{2n+1}$ instead of $T^x_{2n+1}$.

\smallskip

The model ${\bfL}[T_{2n+1}]$ can be seen as an analogue of the
constructible universe ${\bfL}$ in the odd projective levels: The
(Shoenfield) $\Pi^1_1$ scale for a $\Pi^1_1$ complete set is in
${\bfL}$, hence ${\bfL}[T_1] = {\mathbf L}$.\footnote{For a proof,
{\sl cf.} \cite[9C]{KechrisMoschovakis}.} Indeed, the reals of
${\bfL}[T_{2n+1}]$ are exactly the reals of ${\M}_{2n}$, the
canonical iterable inner model with $2n$ Woodin
cardinals.\footnote{Combine \cite[Corollary 4.9]{PWOIM}
with \cite[Theorem 7.2.1]{HarringtonKechris}.}

Moreover, not just the reals of the models, but the models
${\bfL}[T_{2n+1}]$ themselves are independent of the choices of the particular
$\Pi^1_{2n+1}$ complete set and the scale on it,
as has been shown by Becker and Kechris
in \cite[Theorem 1 \& 2]{BeckerKechris}:

\begin{thm}\label{fnbeke} Assume {\sf PD} and let $x\in\reals$ be a real.
If $P$ and $Q$ are $\Pi^1_{2n+1}(x)$ complete sets, $\Phi$ and $\Psi$
are scales on $P$ and $Q$, respectively, and $T$ and $S$ are the
trees associated to $\Phi$ and $\Psi$, respectively. Then $\bfL[T] =
\bfL[S]$.
\end{thm}

\medskip

Another consequence of determinacy which will be mentioned only
briefly to simplify notation is the existence of largest countable
sets of certain (lightface) complexity classes:

\begin{thm}\label{LargestCountable}
Let $x\in\reals$ be a real. Suppose that ${\rm
Det}(\Delta^1_{2n}(x))$ holds. Then there is a largest countable
$\Sigma^1_{2n+2}(x)$ set which will be denoted by
$C_{2n+2}(x)$.\footnote{As for the trees $T_{2n+1}$, we shall omit
the parameter $x$ and write $C_{2n+1}$ if the real $x$ is recursive.}
\end{thm}

\proof {\sl Cf.} \cite[Theorem 2]{KeMo72}.\qed

\section{ON RAMSEY${}^\bigstar$ ULTRAFILTERS}\label{sec:druf}

In this section we show that Ramsey$\sst$ ultrafilters exist if we
assume $\CH$ and that in general both existence and non-existence of
Ramsey$\sst$ ultrafilters are consistent with $\ZFC$.

First we show that an scp-ultrafilter induces in a canonical way a
Ramsey filter on $\omega$.

\begin{fct}\label{fct:MIN}
If $\F\subs\parto$ is an scp-ultrafilter, then
$\{\mMin(X):X\in\F\}\setminus\{0\}$ is a Ramsey filter on
$\omega$.\end{fct}

\proof For positive natural numbers $n$ and $r$ let $\tau:[\omega
]^n\to r$ be any colouring. We define $\pi:\nseg{\omega}{n}\to r$ by
stipulating $\pi (s):=\tau\big(\mMin(s^*)\setminus\{0\}\big)$. It is
easy to see that if $\pi\res{\nseg{X}{n}}$ is constant for an
$X\in\F$, then $\tau\res[\mMin(X)\setminus\{0\}]^n$ is constant, too.
\qed

\begin{prop}
It is consistent with $\ZFC$ that there are no scp-ultrafilters.
\end{prop}

\proof Kunen proved ({\sl cf.}~\cite[Theorem\;91]{Jechbook}) that
it is consistent with $\ZFC$ that there are no Ramsey filters on
$\omega$. Therefore, by Fact~\ref{fct:MIN}, in a model of $\ZFC$ in
which there are no Ramsey filters, there are also no
scp-ultrafilters.\qed

Let $\Us=\la\parto,\le\ra$ be the partial order defined as in Subsection
\ref{subsec:DMF}. It
is easy to see that the forcing notion $\Us$ is $\sigma$-closed (this
is part of Fact~2.3 of~\cite{Lorisym}).

\begin{lm}\label{lm:fURamsey}
If $G$ is $\Us$-generic over ${\V}$, then $G$ is an scp-ultrafilter
in ${\V}[G]$.
\end{lm}

\proof Let $s\in\NN$ and $k\in\omega$ with $|s|=n$ and $n+k>0$.
Further, let $\pi:\nseg{\omega}{(n+k)}\to r$ be any colouring and for
$s\ceq X\in\parto$ let $$H_{\pi (s,X) }
:=\{Y\in\open{s,X}:\pi ({s,X})\res{\nseg{s,Y}{(n+k)}}\text{\ is
constant}\}\,.$$ By the main result of \cite{Loricolour} and its
proof, for every dual Ellentuck neighbourhood $\open{s,X}$ and for
any colouring $\pi: \nseg{s,X}{(n+k)}\to r$, there is a
$Y\in\open{s,X}$ such that $\pi\res\nseg{s,Y}{(n+k)}$ is constant.
Hence, for any dual Ellentuck neighbourhood $\open{s,X}$ and for any
colouring $\pi:\nseg{\omega}{(n+k)}\to r$, the set $H_{\pi (s,X) }$
is dense below $X$. Because every such colouring $\pi$ can be encoded
by a real and $\Us$ is $\sigma$-closed, the forcing notion $\Us$ does
not add any colouring $\pi$, which implies, because $G$ meets each
dense set, that $G$ is an scp-ultrafilter in ${\V}[G]$. \qed

We can prove with similar arguments:

\begin{lm}\label{lm:fUdiag}
If $G$ is $\Us$-generic over ${\V}$, then $G$ is a diagonalizable
ultrafilter in ${\V}[G]$.
\end{lm}

\proof Let $\dot{p}$ be a $\Us$-name such that
$$\force_{\Us}``\dot{p} \;\text{is an $\LU{\bG}$-condition''},$$
where $\bG$ is the canonical name for the $\Us$-generic object, and
let $X$ be any $\Us$-condition. Because $\dot{p}$ can be encoded by a
real number and $\Us$ is $\sigma$-closed, there is a $\Us$-condition
$Y\le X$ and a real $p'\in\V$ such that $Y\force_{\Us} p'=\dot{p}$,
which implies $Y\aceq\suc_{p'}(t)$ for every $t\in p'$. By induction
one can construct a $Z\aceq Y$ such that $Z$ is a branch of $p'$ and
therefore, $$Z\force_{\Us}\text{``there is a branch of $\dot{p}$
which belongs to $\bG$''}.$$ Since $Z\le X$, this completes the
proof. \qed

\begin{prop}\label{prop:CH} Assume $\CH$, then there is a
Ramsey$\sst$ ultrafilter.
\end{prop}

\proof Assume ${\V}\models\CH$. Let $\chi$ be large enough such that
$\cP(\parto)\in H(\chi)$, {\sl i.e.}, the power set of $\parto$ (in
${\V}$) is hereditarily of size $<\chi$. Let $\bfN$ be an elementary
submodel of $\langle H(\chi),\in\rangle$ containing all reals of
${\V}$ with $|\bfN|=2^{\aleph_0}$. We consider the forcing notion
$\Us$ in the model $\bfN$. Because $|\bfN|=2^{\aleph_0}$, in ${\V}$
there is an enumeration $\{D_\alpha\subs\parto:\alpha<2^{\aleph_0}\}$
of all dense sets of $\Us$ which lie in $\bfN$. Since $\Us$ is
$\sigma$-closed and because ${\V}\models\CH$, $\Us$ is
$2^{\aleph_0}$-closed in ${\V}$ and therefore we can construct a
descending sequence $\{p_\alpha:\alpha<2^{\aleph_0}\}$ in ${\V}$ such
that $p_\alpha\in D_\alpha$ for each $\alpha<2^{\aleph_0}$. Let
$G:=\{p\in\parto:p_\alpha\ceq p\ \text{for some $p_\alpha$}\}$, then
$G$ is $\Us$-generic over $\bfN$. By Lemma~\ref{lm:fURamsey} and
Lemma~\ref{lm:fUdiag} we have $\bfN[G]\models\text{``there is a
Ramsey$\sst$ ultrafilter''}$ and because $\bfN$ contains all reals of
${\V}$ and every function $f:\nseg{\omega}{n}\to r$ (where
$n,r\in\omega$) and every $\LU{G}$-condition $p$ can be encoded by a
real number, the Ramsey$\sst$ ultrafilter in $\bfN[G]$ is also a
Ramsey$\sst$ ultrafilter in ${\V}$, which completes the proof. \qed

\section{ON $\LU{\F}$ AND $\MMU{\F}$ FOR RAMSEY${}^\bigstar$
FILTERS $\F$}\label{sec:Laver*}

In this section, $\F\subs\parto$ denotes always a Ramsey$\sst$
ultrafilter.

We shall show that the forcing notions $\LU{\F}$ and $\MMU{\F}$ are
equivalent and that both forcing notions have pure decision and the
homogeneity property (this means that coarsenings of generic objects
remain generic, see Fact~ \ref{HomPureDec}). We show first that
$\MMU{\F}$ has pure decision and the homogeneity property. To show
this we will follow \cite[Section\,4]{Lorisym}.

\medskip

If $s\in\NN$ and $s\ceq X\in\F$, then we call the dual Ellentuck
neighbourhood ${\open{s,X}}$ an {\bf $\F$-dual Ellentuck neighbourhood}
and write ${\opf{s,X}}$ to emphasize that $X\in\F$. A set
$\cO\subs\parto$ is called {\bf $\F$-open} if $\cO$ can be written as
the union of some $\fF$-dual Ellentuck neighbourhoods.

\smallskip

For $s\in\NN$ remember that $s^* =s\cup\big{\{}\{\mdom
(s)\}\big{\}}$.

\smallskip

Let $\cO\subs\parto$ be an $\F$-open set. Call ${\opf{s,X}}$ {\bf
good} (with respect to $\cO$), if for some $Y\in {\opf{s,X}}\cap\F$,
${\opf{s,Y}}\subs\cO$; otherwise call it {\bf bad}. Note that if
${\opf{s,X}}$ is bad and $Y\in {\opf{s,X}}\cap\F$, then ${\opf{s,Y}}$
is bad, too. We call ${\opf{s,X}}$ {\bf ugly} if ${\opf{t^*,X}}$ is
bad for all $s\seg t^*\ceq X$ with $|t|=|s|$. Note that if
${\opf{s,X}}$ is ugly, then ${\opf{s,X}}$ is bad.

\begin{lm}\label{lm:ugly}
Let $\F\subs\parto$ be a Ramsey$\sst$ ultrafilter and
$\cO\subs\parto$ an $\F$-open set. If ${\opf{s,X}}$ is bad (with
respect to $\cO$), then there is a $Z\in {\opf{s,X}}$ such that
${\opf{s,Z}}$ is ugly.
\end{lm}

\proof We begin by constructing an $\LU{\F}$-condition $p$. Let $s_0$
be such that $s\seg s_0^*\ceq X$ and $|s|=|s_0|$, and put
$\stem(p):=s_0$. If there is an $Y\in\opf{s_0^*,X}\cap\F$ such that
$\opf{s_0^*,Y}\subs\cO$, then $X_{s_0}:=Y$, otherwise, $X_{s_0}:=X$.
Let $s_{n+1}^*\seg(s_{n}\kap X_{s_{n}})$ be such that
$|s_{n+1}|=|s_{n}|+1=|s|+n+1$ and let $\{t_i:i\leq h\}$ be an
enumeration of all $t$ such that $s_0\seg t\ceq s_{n+1}$, $|t|=|s|$
and $\mdom (t)=\mdom (s_{n+1})$. Further let $Y^{-1}:=X_{s_n}$. Now
choose for each $i\leq h$ a partition $Y^i\in\F$ such that $Y^i\ceq
Y^{i-1}$, $s_{n+1}^*\seg Y^i$ and ${\opf{{(t_i)^*},Y^i}}$ is bad or
${\opf{{(t_{n}^i)^*},Y^i}}\subs\cO$ and finally, let
$X_{s_{n+1}}:=Y^h$.

\smallskip

Put $p:=\{t\in\NN: s_0\seg t^*\ceq s_n^*$ for some $n\in\omega\}$.
Since $\F$ is diagonalizable, there is a partition $Y\in\F$ which is
a branch of $p$. We may assume that $s_0\seg Y$. Define $S_Y:=\{t:
s\seg t^*\ceq Y\wedge |t|=|s|\}$; then, by the construction of $p$,
for all $t\in S_Y$ we have either ${\opf{t^*,Y}}$ is bad or
${\opf{t^*,Y}} \subs\cO$. Now let $B_0:=\{t\in S_Y: {\opf{t,Y}}$ is
bad$\}$ and $B_1:=\{t\in S_Y:{\opf{t^*,Y}}\subs\cO\}=S_Y\setminus
B_0$. Because $\F$ is an scp-filter, there is a partition $Z\in
{\opf{s,Y}}\cap\F$ such that $S_Z\subs B_0$ or $S_Z\subs B_1$. If we
are in the latter case, we have ${\opf{s,Z}}\subs\cO$, which is a
contradiction to our assumption that ${\opf{s,X}}$ is bad. So, we
must have $S_Z\subs B_0$, which implies that ${\opf{s,Z}}$ is ugly
and completes the proof of the Lemma. \qed

\begin{lm}\label{lm:2ndstep}
If $\F$ is a Ramsey$\sst$ ultrafilter and $\cO\subs\parto$ is an
$\F$-open set, then for every $\F$-dual Ellentuck neighbourhood
${\opf{s,X}}$, there is a $Y\in {\opf{s,X}}\cap\F$ such that
${\opf{s,Y}}\subs\cO$ or ${\opf{s,Y}}\cap\cO\cap\fF =\emptyset$.
\end{lm}

\proof If ${\opf{s,X}}$ is good, then we are done. Otherwise, we can
construct an $\LU{\F}$-condition $p$ in a similar way as in
Lemma~\ref{lm:ugly}, such that for any branch $Y$ of $p$ which
belongs to $\F$ we have the following: For each $t$ with $s\seg
t^*\ceq Y$, the set $\opf{t^*,Y}$ is bad. We claim that
${\opf{s,Y}}\cap\cO\cap\F =\emptyset$. Take any $Z\in
{\opf{s,Y}}\cap\cO\cap\fF$. Because $\cO$ is $\F$-open we find a
$t\seg Z$ such that ${\opf{t^*,Z}}\subs\cO$. Because $s\seg t^*\ceq
Y$ we have by construction that ${\opf{t^*,Y}}$ is bad. Hence, there
is no $Z\in {\opf{t^*,Y}}$ such that ${\opf{t^*,Z}}\subs\cO$. This
completes the proof. \qed

Now we can show that $\MMU{\F}$ has pure decision and the homogeneity
property.

\begin{thm}\label{thm:pd}
Let $\fF$ be a Ramsey$\sst$ ultrafilter and let $\Phi$ be a sentence
of the forcing language $\MMU{\F}$. For any $\MMU{\F}$-condition $\la
s,X\ra$ there is a $\MMU{\F}$-condition $\la s,Y\ra\le\la s,X\ra$
such that $\la s,Y\ra\force_{\MMU{\F}} \Phi$ or $\la
s,Y\ra\force_{\MMU{\F}} \neg\Phi$.
\end{thm}

\proof The proof is same as the proof of
\cite[Theorem\,4.3]{Lorisym}, using Lemma \ref{lm:2ndstep}.\qed

The next theorem shows in fact that if $\F$ is a Ramsey$\sst$
ultrafilter, then $\MMU{\F}$ is proper.

\begin{thm}\label{thm:proper}
Let $\F\subs\parto$ be a Ramsey$\sst$ ultrafilter, then $\MMU{\F}$
has the homogeneity property.
\end{thm}

\proof The proof is same as the proof of
\cite[Theorem\,4.4]{Lorisym}, using Lemma \ref{lm:2ndstep}.\qed

In order to show that $\MMU{\F}$ and $\LU{\F}$ are equivalent if $\F$
is a Ramsey$\sst$ ultrafilter, we define first some special
$\LU{\F}$-conditions.

An $\LU{\F}$-condition $p$ is called {\bf uniform} if there is a
partition $X\in\F$ such that $\open{t,X}=\open{t,X_t^p}$ for every
$t\in p$; this partition is denoted by $\cu (p)$. These conditions
roughly correspond to
the simple conditions of \cite[Definition 1.10]{JudahShelah.12}.

\begin{lm}\label{lm:dense}
If $\F$ is a Ramsey$\sst$ ultrafilter, then the set of all uniform
$\LU{\F}$-conditions is dense and open in $\lU{\F}$.
\end{lm}

\proof Let $p\in\lU{\F}$ with $s=\stem(p)$, then, since $\F$ is
diagonalizable, there is an $X\in\F$ which is a branch of $p$. Let
$q$ be the uniform condition with $\cu (q)=X$ and $\stem(q)=s$. Note
that $X$ is a branch of $q$. By Fact~\ref{fct:sub-branch}, each
$Y\in\open{s,X}$ is also a branch of $p$, which implies that $q\le
p$. \qed

\begin{thm}
If $\F$ is a Ramsey$\sst$ ultrafilter, then $\MMU{\F}$ and $\LU{\F}$
are forcing equivalent.
\end{thm}

\proof Let $I:=\{p\in\lU{\F} :p\ \text{is uniform}\}$ and define
$$\begin{array}{rccc} j: & I & \longrightarrow & M^{\st}_{\F} \\
   & p & \longmapsto     & \la \stem(p),\cu (p)\ra\,,
\end{array}$$
then it is easily checked that $j$ is a dense embedding
and because (by
Lemma~\ref{lm:dense}) $I$ is dense open in $\lU{\F}$, this completes
the proof. \qed

This is the promised dualization of Theorem \ref{LavMath} and possibly one step
towards a proof of Conjecture \ref{OurConjecture}.

\section{THE DUAL RAMSEY PROPERTY FOR SIMPLE POINTCLASSES}
\label{sec:simple}

In the following we will show that it is consistent with $\ZFC$ that the
sets in the first levels of the projective hierarchy are dual
Ramsey. We begin with the analytic sets:

Because $\MMs$ has pure decision and the homogeneity property, one
can show the pretty straightforward

\begin{fct}\label{Sigma-1-1-Dual-Ramsey}
Analytic sets are completely dual Ramsey.
\end{fct}

\proof Let $A$ be an arbitrary $\Sigma^1_1(a)$ set with parameter
$a\in\reals$ and let $\open{s,Y}$ be any dual Ellentuck neighbourhood
and $\la s,Y\ra$ the corresponding $\MMs$-condition. Take a countable
model $\bfN$ of a sufficiently large fragment of $\ZFC$ which contains
$Y$ and $a$. Let $\Gen$ be the canonical name for the $\MMs$-generic
object. Because $\MMs$ has pure decision we find an $\MMs$-condition
$\la s,Z\ra\le\la s,Y\ra$ which decides ``$\Gen\in \dot{A}$''. Since
$\bfN$ is countable, there is an $X\in\open{s,Z}$ which is
$\MMs$-generic over $\bfN$ and because every $X'\in\open{s,X}$ is also
$\MMs$-generic we have $$\bfN[X]\models\text{``$\open{s,X}\subs A$ or
$\open{s,X}\cap A=\emptyset$''}.$$ Because $A$ and $\open{s,Y}$ were
arbitrary and $\BSigma^1_1$ sets are absolute between $\V$ and $\bfN$,
we are done. \qed

Note that Fact~\ref{Sigma-1-1-Dual-Ramsey} is verified without any
reference to forcing by looking at the unfolded version of the
Banach--Mazur game for the dual Ellentuck topology.\footnote{{\sl
Cf.} \cite[Theorem (21.8)]{Kechris}.}

\bigskip

Remember that according to Proposition \ref{prop:Delta.1n}
if every $\BDelta^1_2$ set is dual Ramsey then every
$\BDelta^1_2$ set has the classical Ramsey property. Because it is
not provable in $\ZFC$ that every $\BDelta^1_2$ set is Ramsey, it is
also not provable in $\ZFC$ that every $\BDelta^1_2$ set is dual
Ramsey ({\sl e.g.}, $\bfL\models$ ``There is a $\BDelta^1_2$ set
which is not dual Ramsey''). On the other hand we have

\begin{fct} The following theories are equiconsistent:
\begin{itemize}
\item[(a)] {\ZFC}.
\item[(b)] {\ZFC}$~+~\CH~+$ every $\BSigma^1_2$ set is dual Ramsey.
\item[(c)] {\ZFC}$~+~2^{\aleph_0}=\aleph_2~+$ every $\BSigma^1_2$
set is dual Ramsey.
\end{itemize}
\end{fct}

\proof We get both (b) and (c) by an iteration (of length $\omega_1$
and $\omega_2$ respectively) of dual Mathias forcing with countable
support, starting from $\bfL$ ({\sl cf.} also
\cite[Theorems 6.2 \& 6.3]{Lorisym}). \qed

Remember that $\fJ=\{A\subs\parto: A\ \text{is completely dual
Ramsey-null}\}$ and let

$$\add\fJ:=\mmin\{|\fE |: \fE\subs\fJ\wedge\bigcup\fE\notin\fJ\}$$

and

$$\cov\fJ:=\mmin\{|\fE |: \fE\subs\fJ\wedge\bigcup\fE
=\parto\}\,.$$

In \cite{Lorisha} it is shown that $\add\fJ=\cov\fJ=\fH$, where $\fH$
is the dual shattering cardinal. If $\cJ$ denotes the ideal of
classical completely Ramsey-null sets, then we get the analogous
result, namely $\add\cJ=\cov\cJ=\fh$, where $\fh$ is the shattering
cardinal ({\sl cf.}~\cite{Plewik}). Because every $\BSigma^1_2$ set
is the union of $\aleph_1$ Borel sets ({\sl
cf.}~\cite[Theorem\,95]{Jechbook}), it is easy to see that $\fH
>\aleph_1$ implies that every $\BSigma^1_2$ set is even completely
dual Ramsey (and the analogous result holds for the classical Ramsey
property with respect to $\fh$). Now, an $\omega_2$-iteration with
countable support of dual Mathias forcing starting from $\bfL$ yields
a model in which $\fH=\aleph_2$ ({\sl cf.}~\cite{Lorisha}). Thus,
this provides another proof that ``Every $\BSigma^1_2$ set is dual
Ramsey'' is consistent with $\ZFC$. In Section~\ref{sec:det} we shall
provide a third proof as a byproduct of the analysis of scales under
{\sf PD}.

\medskip

Concerning Martin's Axiom $\MA$, it is well-known that $\MA$ implies
$\fh=2^{\aleph_0}$. Hence, by the facts mentioned above, $\MA +
\neg\CH$ implies that all $\BSigma^1_2$ sets have the classical
Ramsey property.

\medskip

A similar argument for the dualized case does not work: Brendle has
shown in \cite{BrendleMA}, that $\MA + 2^{\aleph_0}>\fH=\aleph_1$ is
consistent with $\ZFC$. However, this result does not preclude that
$\MA +\neg\CH$ might imply that every $\BSigma^1_2$ set is dual Ramsey.

\bigskip

At the next level we like to mention the following

\begin{fct} ``Every $\BDelta^1_3$ set is dual Ramsey'' is consistent
with $\ZFC$.
\end{fct}

\proof An $\omega_1$-iteration (or also an $\omega_2$-iteration) with
countable support of dual Mathias forcing starting from $\bfL$ yields
a model in which every $\BDelta^1_3$ set is dual Ramsey. The proof is
exactly the same as the proof of the corresponding result for the
classical Ramsey property given in \cite{JudahShelah.13}; the reason
for this is that they needed only that Mathias forcing is proper and
has the homogeneity property, but these two properties hold also for
dual Mathias forcing. \qed

\bigskip

The $\boldsymbol{\Delta}^1_3$ level is probably as far as we can get
in {\ZFC} without further assumptions. It is a famous open question
whether the consistency of ``Every $\boldsymbol{\Pi}^1_3$ set is
Ramsey'' implies the existence of an inner model with an inaccessible
cardinal ({\sl cf.} \cite[Question 11.16]{Kan94} \&
\cite[p.~49]{Raisonnier}). Most likely, the dualized question is
equally hard to conquer.

\section{DETERMINACY AND THE DUAL RAMSEY PROPERTY}\label{sec:det}

We shall move on to arbitrary projective sets in this section. As we
mentioned earlier, this means that we probably have to go beyond
\ZFC.

In \cite[Section\;5]{CarlsonSimpson}, the authors prove in fact that
in the Solovay model constructed by collapsing an inaccessible
cardinal to $\omega_1$ every projective set is dual Ramsey. As we remarked,
it is unknown whether the inaccessible
cardinal is necessary for that.

\medskip

But there is another question connected to the dual Ramsey property
of projective sets: As with the standard Ramsey property we can ask
whether an appropriate amount of determinacy implies the dual Ramsey
property. As usually with regularity properties of sets of reals we
would expect that ${\rm Det}(\BPi^1_n)$ implies the dual Ramsey
property for all $\BSigma^1_{n+1}$ sets. But a direct implication
using determinacy is not as easy as with the more prominent
regularity properties (as Lebesgue measurability and the Baire
property) since the games connected to the dual Ramsey property (the
Banach--Mazur games in the dual Ellentuck topology) cannot be played
using natural numbers.

\smallskip

The same problem had been encountered with the standard Ramsey
property and had been solved in \cite{HarringtonKechris} by making use of
the scale property and the Periodicity Theorems \ref{FPT} and \ref{SPT}:

\begin{thm}
If ${\rm Det}(\BDelta^1_{2n+2})$, then every $\BPi^1_{2n+2}$ set is Ramsey.
\end{thm}

The main ingredient of this proof was an analysis of the models
${\bfL}[T_{2n+1}]$
under Determinacy assumptions (Lemma
\ref{LTsmall}). In the following we shall give a brief review of the
result with sketches of
an adaptation to our context.

\bigskip

\begin{lm}\label{gettrees}
Let $\Gamma$ be any $\omega$--parametrized pointclass. If $U$ is
$\omega$--universal for $\Gamma$, $A\in\Gamma$, $\bfN$ any model, and
$T\in\bfN$ a tree such that $p[T] = U$. Then there is a tree
$S\in\bfN$ such that $p[S] = A$.
\end{lm}

\proof This is basically \cite[Proposition 13.13 (g)]{Kan94}, apart
{}from the assertion that $S\in\bfN$. But this is clear since the
reduction function reducing $A$ to $U$ is just the trivial function
$x\mapsto\la n_0,x\ra$ (where $n_0$ is the index of $A$ in $U$) and
hence in $\bfN$.\qed

\begin{lm}\label{LTsmall}
Assume ${\rm Det}(\BDelta^1_{2n+2})$. Then
$\reals\cap{\bfL}[T_{2n+1}] = C_{2n+2}$. In particular, this is a
countable set.
\end{lm}

{\it Proof Sketch.} We shall very roughly sketch the argument of
\cite[Theorem 7.2.1]{HarringtonKechris}:

First of all $\reals\cap\bfL[T_{2n+1}]$ is easily seen to be
$\Sigma^1_{2n+2}$. That every countable $\Sigma^1_{2n+2}$ set of
reals is a member of $\bfL[T_{2n+1}]$ follows directly from
Mansfield's Theorem ({\sl cf.} \cite[Theorem 14.7]{Kan94}) and Lemma
\ref{gettrees}.\footnote{Note that by Theorem
\ref{fnbeke} the choice of the complete
set for the definition of the $\bfL[T_{2n+1}]$ doesn't matter.}

So, what is left to show is that $\reals\cap\bfL[T_{2n+1}]$ actually
is countable. The proof uses the following steps:

\begin{enumerate}

\item Fix a $\Pi^1_{2n+1}$ norm $\varphi^* : P_{2n+1} \to
\boldsymbol{\delta}^1_{2n+1}$ which exists according to Theorem
\ref{FPT}.

\item Using $\phi^*$, code the tree $T_{2n+1}$ by some $A\subseteq
\boldsymbol{\delta}^1_{2n+1}$ and show that ${\bfL}[T_{2n+1}] = {\bfL}[A]$.

\item For arbitrary subsets $X\subseteq\boldsymbol{\delta}^1_{2n+1}$, define
$X^* := \{ z\in\reals~:~\phi^*(z) \in X\}$.

\item Show: If $X^*\in\Delta^1_{2n+2}$, then the set $\reals\cap{\bfL}[X]$ is
contained in a countable $\Sigma^1_{2n+2}$ set.

\item Compute: $A^*\in\Delta^1_{2n+2}$.
\end{enumerate}\qed

Harrington and Kechris used this result to receive
results about projective sets from {\sf PD} alone that formerly could
only be derived from stronger hypotheses. The results for the classical Ramsey
property follows Solovay's argument for the $\BSigma^1_2$ case.
We shall outline this
argument in full generality and then apply it to the dual Ramsey
property.

\medskip

At first we need to relativize Lemma \ref{LTsmall} in two different parameters:

\begin{lm}\label{LTsmallrel}
Assume ${\rm Det}(\BDelta^1_{2n+2})$. Let $x,y\in\reals$ be any real numbers.
Then $\reals\cap\bfL[T_{2n+1}^y,x] = C_{2n+2}(x\oplus y)$.
\end{lm}

{\it Proof Sketch.}
As an immediate relativization of
of Lemma
\ref{LTsmall} (for the pointclass $\Pi^1_{2n+1}(y)$ instead of $\Pi^1_{2n+1}$),
we get:

$$\reals\cap\bfL[T^y_{2n+1}] = C_{2n+2}(y).$$

\smallskip

To show that $\reals\cap\bfL[T^y_{2n+1},x]$ is a countable set, we
have to relativize (iv) and (v) again. The obvious relativization of
(iv) is

\begin{itemize}
\item[(iv*)] If $X^*\in\Delta^1_{2n+2}(x\oplus y)$,
then the set $\reals\cap{\bfL}[X]$ is
contained in a countable $\Sigma^1_{2n+2}(x\oplus y)$ set.
\end{itemize}

Since $\bfL[T^y_{2n+1}] = \bfL[A]$ for some $A\subseteq
\boldsymbol{\delta}^1_{2n+1}$ such that $A^*$ is $\Delta^1_{2n+2}(y)$
according to (ii.) and (v.), we know that $\bfL[T^y_{2n+1},x] =
\bfL[A,x]$. Thus we have to find a set $B\subseteq
\boldsymbol{\delta}^1_{2n+1}$ such that $\bfL[B] = \bfL[A,x]$ and
$B^*\in\Delta^1_{2n+1}(x\oplus y)$. This would prove the theorem.

\medskip

The natural choice for $B$ is:
\begin{eqnarray*}
B & := & \big\{\alpha \in\boldsymbol{\delta}^1_{2n+1} :
\big(\alpha\geq\omega\wedge\alpha\in A\big) \mbox{ or } \\ &&
\big(\alpha<\omega\wedge\exists n (\alpha = 2n\wedge n\in
A)\big)\mbox{ or }\\ && \big(\alpha<\omega\wedge\exists n (\alpha =
2n+1\wedge n\in x)\big)\big\}.
\end{eqnarray*}

Obviously, $B\in\bfL[A,x]$ and $A$ and $x\in\bfL[B]$, so $\bfL[A,x]=\bfL[B]$.
Thus, what is left is to show
that $B^*$ is $\Delta^1_{2n+2}(x\oplus y)$.
But this is easy to see using the definition of
$B$, and the facts that $\phi^*$ was a $\Pi^1_{2n+1}$ norm and that $A^*$ was
$\Delta^1_{2n+2}(y)$.\qed

\medskip

The following lemma is an obvious generalization of Shoenfield's Absoluteness
Lemma ({\sl cf.} also \cite[Theorem 8G.10]{Moschovakis}):

\begin{lm}\label{TAbs}
$\Sigma^1_{2n+2}(x)$ formulae are absolute for models containing
$T^x_{2n+1}$, {\sl i.e.}, if $\bfN$ is a model with
$T^x_{2n+1}\in\bfN$ and $\phi$ is a $\Sigma^1_{2n+2}$ formula, then
$$\forall X\in\bfN~(\bfN\models \phi[X,x]~\iff~\V\models\phi[X,x]).$$
\end{lm}

\proof By
Theorem \ref{fnbeke}, we can assume
that $T^x_{2n+1}$ was constructed using an $\omega$--universal set for
$\Pi^1_{2n+1}(x)$, enabling us to use Lemma \ref{gettrees}.

\smallskip

Thus every $\Pi^1_{2n+1}(x)$ set is represented by a tree $S\in\bfN$.
We easily get a tree $S^*$ for each $\Sigma^1_{2n+2}(x)$ set ({\sl
cf.} \cite[Proposition 13.13 (d)]{Kan94}).

But now the theorem follows from standard absoluteness of illfoun\-ded\-ness
as in Shoenfield's proof ({\sl cf.}
\cite[Exercise 12.9 (a)]{Kan94}).\qed

\begin{thm}\label{Ramsey}
Let $x\in\reals$ be a real.
Suppose that there is a dual Mathias generic partition
over $\bfL[T^x_{2n+1}]$. Then every $\Sigma^1_{2n+2}(x)$ set is dual Ramsey.
\end{thm}

\proof Let $A$ be a $\Sigma^1_{2n+2}(x)$ set and $\phi$ a $\Sigma^1_{2n+2}(x)$
expression
describing $A$, {\sl i.e.}
$$\forall y~(~y\in A ~\iff~ \phi[y]~).$$

\smallskip

By Fact \ref{HomPureDec} we find a $\MMs$--condition $\la
\emptyset,X\ra\in\bfL[T^x_{2n+1}]$
such that
$$\mbox{either }\la \emptyset,X\ra\force \phi(\pc(\bXG))\mbox{ or }\la
\emptyset,X\ra\force
\neg\phi(\pc(\bXG)),$$ where $\bXG$ is the name for a dual Mathias
 generic partition.

\smallskip

Without loss of generality, we assume the former. By our assumption,
we actually have a generic partition $Z$ over $\bfL[T^x_{2n+1}]$ with
$\la \emptyset,X\ra\in G_Z$, where $G_Z$ is the filter associated to
$Z$, {\sl i.e.}

$$\la t,Y\ra\in G_Z ~\iff~ Z\in\open{t,Y}.$$

This means that $Z\in\open{X}$. Now by \ref{HomPureDec} (homogeneity
of $\MMs$) again, every element $Z^*$ of $\open{Z}$ is also
$\MMs$--generic over $\bfL[T_{2n+1}^x]$. Since $G_{Z^*}\subseteq
G_Z$, we still have $\la \emptyset,X\ra\in G_{Z^*}$. Consequently, we
have $\bfL[T^x_{2n+1}][Z^*]\models \phi[\pc(Z^*)]$. But $\phi$ was
absolute for models containing $T^x_{2n+1}$ by Lemma \ref{TAbs},
hence we have $\V\models\phi[\pc(Z^*)]$.

\smallskip

Summing up, we have found a partition $Z$ such that $\{\pc(Z^*) :
Z^*\ceq Z\}\subseteq A$. This is exactly what we had to show by
Observation \ref{INVERSE}.\qed

Note that this type of argument probably will not work if you replace
``dual Ramsey'' by ``completely dual Ramsey''.
What you would have to do is to relativize the argument to arbitrary partitions
$W\in\V$. But
at least this does not work
in the classical case: Brendle has shown in \cite{BrendlePrague} that
in any model containing one Mathias real over a ground model $\bfN$,
there is an Ellentuck neighbourhood that doesn't contain any Mathias
reals over $\bfN$.

\medskip

Another useful comment about Theorem \ref{Ramsey} is that if you look
at the case $n=0$ you get a third proof of the consistency of ``Every
$\BSigma^1_2$ set is dual Ramsey'':

\begin{cor}\label{GenericsSigma12}
Suppose that for each real $x\in\reals$ there is a dual Mathias
generic partition over $\bfL[x]$. Then every $\BSigma^1_2$ set is
dual Ramsey.
\end{cor}

\proof Immediate from Theorem \ref{Ramsey}, keeping in mind that
$T_1\in\bfL$ by \cite[9C]{KechrisMoschovakis}, as mentioned in
Subsection \ref{subsec:cabal}.\qed

This particularly generic version of proving the consistency of
properties of $\BSigma^1_2$ sets should be compared to analogous
results for Random forcing, Cohen forcing and Hechler
forcing.\footnote{{\sl Cf.} \cite{BrLoe99}. Note that in most cases
the existence of generics doesn't give more than regularity at the
$\BDelta^1_2$ level, and something more than mere existence
is needed for the $\BSigma^1_2$ level.}

\bigskip

We now move on to use Lemma \ref{LTsmallrel} and Theorem \ref{Ramsey} to
get that Determinacy implies the dual Ramsey property:

\begin{cor} \label{PDDualRamsey}
Assume ${\rm Det}(\BDelta^1_{2n+2})$.
Then every $\BSigma^1_{2n+2}$ set is dual Ramsey.
\end{cor}

\proof By Lemma \ref{LTsmallrel}, we get that for any reals
$x,y\in\reals$, the set of reals in $\bfL[T^x_{2n+1},y]$ is
countable.

By the same argument that is used to show that $\forall x
(\aleph_1^{\bfL[x]} < \aleph^\V_1)$ implies that $\aleph_1^\V$ is
strongly inaccessible in every $\bfL[x]$, we get that
$\cP(\parto)\cap\bfL[T_{2n+1}^x,y]$ is countable for arbitrary
choices of $x$ and $y\in\reals$.

Thus there are dual Mathias generic partitions over each
$\bfL[T^x_{2n+1},y]$, in particular over each $\bfL[T^x_{2n+1}]$, and
we can use Theorem \ref{Ramsey} to prove the claim.\qed

\section{Appendix: Game-filters have the segment-colouring-property}
\label{sec:appendix}

Let $\F\subs\parto$ be an ultrafilter. Associated with $\F$ we define
the game $\cG_{\F}$ as follows. This type of game, which is the
Choquet-game with respect to the dual Ellentuck topology ({\sl
cf.}~\cite[8.C]{Kechris}), was first suggested by Kastanas in
\cite{Kastanas}.

$$\begin{array}{ccccccccc} {\operatorname{I}} &\ \ \ &\la t_0,Y_0\ra
& &\la t_1,Y_1\ra & &\la t_2,Y_2\ra & &  \\
   &      &            &   &           &   & &          & \ldots \\
{\operatorname{II}}& & &\la X_0\ra & &\la X_1\ra & &\la X_2\ra &
\end{array}$$

All the moves $X_n$ of player~II must be elements of the ultrafilter
$\F$ and all the moves $\la t_n,Y_n\ra$ of player~I plays must be
such that $Y_n\in\fF$ and $\open{t_n^*,Y_n}$ is a dual Ellentuck
neighbourhood. Further, the $n$th move $X_n$ of player~II is such
that $X_n\in\open{t_{n}^*,Y_{n}}$ and then player~I plays $t_{n+1}$
such that $t_{n}^*\seg t_{n+1}^*\ceq X_n$ and
$|t_{n+1}|=|t_{n}|+1=|t_0|+n+1$. Player~I wins if and only if the
unique $Y$ with $t_n\seg Y$ (for all $n$) is not in $\F$.

\smallskip

An ultrafilter $\F$ is a {\bf game filter} if and only if player~I has no
winning strategy in the game $\cG_{\F}$.\footnote{For
the existence of game filters see \cite{Lorisym}, where one
can find also some results concerning dual Mathias forcing restricted
to such filters.}

\medskip

In the following we outline the proof that game filters are also
scp-filters. The crucial point will be to show the Preliminary Lemma
\ref{PreliminaryLemma}, which is in fact Carlson's Lemma ({\sl
cf.}~\cite[Lemma~2.4]{CarlsonSimpson}) restricted to game filters.
But first we have to give some notations.

Let $s,t\in\NN$ be such that $s\ceq t$, $|s|=n$ and $|t|=m$. For $k$
with $k\le m-n$ let
$$\nopeni{t}{k}{s}:=\{u\in\NN:\mdom(u)=\mdom(t)\wedge s\seg u\ceq
t\wedge |u|=|s|+k\}\,.$$ For $s\seg t\ceq X$, let
$$\nopeni{t,X}{k*}{s}:=\{u\in\NN: t\seg u^*\ceq X\wedge
|u|=|s|+k\}\,,$$ and let $\nopeni{X}{k*}{s}:=\nseg{s,X}{(n+k)}$.

We have chosen this notation following
\cite[Definition\,2.1]{CarlsonSimpson}, where one can consider $s$ as
an alphabet of cardinality $n$.

For the remainder of this section, let $\F$ be an arbitrary but fixed
game filter.

\begin{plm} \label{PreliminaryLemma}
Let $s\seg X\in\F$ and $\pi:\nopeni{X}{0*}{s}\to l$, then there
exists a $Y\in\open{s,X}\cap\F$ such that
$\pi\res{\nopeni{Y}{0*}{s}}$ is constant.
\end{plm}

Following the ideas of the proof of Theorem~6.3 of
\cite{CarlsonSimpson}, the proof of the Preliminary Lemma will be
given in a sequence of lemmas. We start by stating the well-known
Hales--Jewett Theorem in our notation.

\begin{HJthm} Let $s\in\NN$. For all $d\in\omega$, there is an
$h\in\omega$ such that for any $t\in\NN$ with $s\seg t$ and
$|t|=|s|+h$, and for any colouring $\tau:\nopeni{t}{0}{s}\to d$,
there is a $u\in\nopeni{t}{1}{s}$ such that ${\nopeni{u}{0}{s}}$ is
monochromatic.
\end{HJthm}

The number $h$ in the Hales--Jewett Theorem depends only on the
number $d$ and the size of $|s|$. Let HJ$(d,|s|)$ denote the smallest
number $h$ which verifies the Hales--Jewett Theorem.

Let $s,t\in\NN$ and $X\in\F$ be such that $s\seg t\ceq X$. A set
$K\subs\NN$ is called {\bf dense in} $\open{t,X}$, if for all
$Y\in\open{t,X}\cap\F$, there is a $u$ with $t\seg u^*\ceq Y$ which
belongs to $K$. A set $D\subs\nopeni{\omega}{k*}{s}$ is called {\bf
$k$-dense in} $\nopeni{t,X}{k*}{s}$, if for all
$Y\in\open{t,X}\cap\F$, we have $\nopeni{t,Y}{k*}{s}\cap D\neq
\emptyset$.

\begin{lm}\label{lm:lclaim}
Let $s\seg t\seg X\in\F$ and assume that
$D\subs\nopeni{\omega}{0*}{s}$ is $0$-dense in $\nopeni{X}{0*}{s}$.
Further assume that $K=\{u: t\seg u \wedge \nopeni{u}{0}{s}\cap
D\neq\emptyset\}$ is dense in some $\open{t,Z}$, where
$Z\in\open{t,X}\cap\F$. Then there is an ${\bar s}^*\ceq Z$ with
$t\seg\bar s$ such that for all $v^*\ceq Z$ with ${\bar s}\seg v$ we
have $\nopeni{v}{0}{s}\cap D\neq\emptyset$.
\end{lm}

\proof We shall define a strategy for player~I in the game $\cG_{\F}$,
such that player~I can follow this strategy just in the case when
Lemma~\ref{lm:lclaim} fails. This means that for every ${\bar
s}^*\ceq Z$ with $t\seg\bar s$ there is a $v^*\ceq Z$ with ${\bar
s}\seg v$ such that $\nopeni{v}{0}{s}\cap D= \emptyset$.

Let $t_0^*\ceq Z$ be such that $t\seg t_0$, $|t|=|t_0|$ and
$\nopeni{t_0}{0}{s}\cap D=\emptyset$. Further put $Y_0=t_0\kap Z$ and
player~I plays $\la t_0,Y_0\ra$. Assume $\la t_m,Y_m\ra$ is the $m$th
move of player~I and player~II replies with $\la X_{m}\ra$. If the
lemma fails with $\bar s =t_{m}^*$, player~I can play $\la
t_{m+1},Y_{m+1}\ra$, according to the rules of the game, such that
$\nopeni{t_{m+1}}{0}{s}\cap D=\emptyset$.

\smallskip

Since $\F$ is a game filter, the strategy of player~I is not a
winning strategy and the unique $Y\in\open{Z}$ such that $t_m\seg Y$
(for all $m\in\omega$) belongs to $\F$. Take an arbitrary $u$ with
$t\seg u^*\ceq Y$. For such a $u$ we find a $t_n\seg Y$ such that
$u\ceq t_n$ and $\mdom (u)= \mdom (t_n)$. By the strategy of player~I
we have $\nopeni{t_n}{0}{s}\cap D=\emptyset$ and therefore
$\nopeni{u}{0}{s}\cap D=\emptyset$. But this is a contradiction to
the assumption that $K$ is dense in $\open{t,Z}$. Hence, player~I
cannot follow this strategy, which completes the proof. \qed

\begin{lm}\label{lm:l65}
Suppose $s\seg t\seg X\in\F$, $D$ is $0$-dense in $\nopeni{X}{0*}{s}$
and $K=\{u: t\seg u\;\wedge\;\nopeni{u}{0}{s}\cap D\neq\emptyset\}$
is dense in $\open{t,Z}$, where $Z\in\open{t,X}\cap\F$. Then there is
a ${\bar t}^*\ceq Z$ with $t\seg \bar t$ and $|\bar t|=|t|+1$ such
that $\nopeni{\bar t}{0}{s}\subs D$.
\end{lm}

\proof Let $\bar s$ be as in the Lemma~\ref{lm:lclaim}, and let
$d:=|\nopeni{\bar s}{0}{s}|$. By the Hales--Jewett Theorem, let
$h:=\operatorname{HJ}(d,|s|)$. Pick $v\in\NN$ such that $\bar s\seg
v^*\ceq X$ and $|v|=|\bar s|+h$. Let $\{s_i:s_i\in\nopeni{\bar
s}{0}{s}\wedge i\in d\}$ be an enumeration of the elements of
$\nopeni{\bar s}{0}{s}$. We colour $\nopeni{v}{0}{s}$ by stipulating
$\tau (u)=i$ if and only if $s_i\seg u\wedge u\in D$. By the choice
of $n$, there are $\bar t\in\nopeni{v}{1}{s}$ such that $\nopeni{\bar
t}{0}{s}$ is monochromatic, and therefore, $\nopeni{\bar
t}{0}{s}\subs D$. Thus, we have found a $\bar t$ with $t\seg\bar t$
and $|\bar t|=|t|+1$ such that $\nopeni{\bar t}{0}{s}\subs D$. \qed

\begin{lm}\label{lm:l66}
Suppose $s\seg X\in\F$ and $D$ is $0$-dense in $\nopeni{X}{0*}{s}$.
Then there are $t\in\nopeni{X}{1*}{s}$ and $Y\in\open{t,X}\cap\F$
such that $\{u:t\seg u\wedge\nopeni{u}{0}{s}\subs D\}$ is $1$-dense
in $\nopeni{t,Y}{1*}{s}$.
\end{lm}

\proof In a similar way as above we can define a strategy for
player~I in the game $\cG_{\F}$, such that player~I can follow this
strategy only if Lemma~\ref{lm:l65} fails. But if Lemma~\ref{lm:l65}
is wrong, this would yield\;--\;because $\F$ is a game filter\;--\;a
contradiction ({\sl cf.}~\cite[Lemma\,6.5]{CarlsonSimpson}). \qed

Notice that in Lemma~\ref{lm:l66} we did not require that the
$r\in\NN$ for which we have $r^*\seg t\in\nopeni{X}{1*}{s}$ belongs
to $D$. This we do in

\begin{lm}\label{lm:l67}
Suppose $s\seg X\in\F$ and $D$ is $0$-dense in $\nopeni{X}{0*}{s}$.
Then there are $r^*\seg t\in\nopeni{X}{1*}{s}$ and
$Y\in\open{t,X}\cap\F$ such that $\{u:t\seg
u\wedge\nopeni{u}{0}{s}\subs D\}$ is $1$-dense in
$\nopeni{t,Y}{1*}{s}$ and $r\in D$.
\end{lm}

\proof Let $Y_0\in\open{t_0,X}\cap\F$ be as in the conclusion of
Lemma~\ref{lm:l66}. Thus $D_0:=\{u:t_0\seg u\wedge
\nopeni{u}{0}{s}\subs D\}$ is $1$-dense in $\nopeni{t_0,Y_0}{1*}{s}$.
Using Lemma~\ref{lm:l66}, player~I can play  $\la t_m,Y_m\ra$ at the
$m$th move such that $D_m:=\{u:t_m\seg u\wedge \nopeni{u}{0}{s}\subs
D\}$ is $(m+1)$-dense in $\nopeni{t_m,Y_m}{(m+1)*}{s}$.

Because player~I has no winning strategy, the unique $Y\in\open{s,X}$
such that $t_m\seg Y$ (for all $m$) belongs to $\F$, and because $D$
is $0$-dense in $\nopeni{X}{0*}{s}$, there is an
$r\in\nopeni{Y}{0*}{s}$ which belongs to $D$. Let $\bar
t\in\nopeni{t_m}{1}{s}$ be such that $r^*\seg \bar t$. Since
$\nopeni{\bar t}{0}{s}\subs\nopeni{t_m}{0}{s}$ and because $D_m$ is
$(m+1)$-dense in $\nopeni{t_m,Y_m}{(m+1)*}{s}$ we get $\{u:\bar t\seg
u\wedge\nopeni{u}{0}{s}\subs D\}$ is $1$-dense in
$\nopeni{\bar{t},Y_m}{1*}{s}$. Hence, we have found an $r^*$ such
that $\{u:r^*\seg u\wedge\nopeni{u}{0}{s}\subs D\}$ is $1$-dense in
$\nopeni{r^*,Y}{1*}{s}$ and $r\in D$. \qed

\bigskip

Now we can go back to the

\smallskip

\poplm Let $s\seg X\in\F$. We have to show that for any colouring
$\pi:\nopeni{X}{0*}{s}\to l$, there is a $Y\in\open{s,X}\cap\F$
such that $\nopeni{Y}{0*}{s}$ is monochromatic.

It is easy to see that at least one of the colours is $0$-dense in
$\nopeni{X}{0*}{s}$, say $j$ and let $D:=\{t\in
\nopeni{X}{0*}{s}:\pi(t)=j\}$. Now we can prove the Preliminary~Lemma
in almost the same way as Lemma~\ref{lm:l67}, the only difference is
that player~I uses now Lemma~\ref{lm:l67} to construct the $m$th
move, instead of Lemma~\ref{lm:l66}. \qed

Finally we get the main result of this section.

\begin{prop}\label{prop:games}
Each game filter is also an scp-filter.
\end{prop}

\proof We have to show that for any colouring
$\pi:\nseg{s,\omega}{(|s|+k)}\to r$, where $r$ and $k$ are positive
natural numbers and $s\in\NN$, there is an $X\in\F$ such that
$s\seg X$ and $\nseg{s,X}{(|s|+k)}$ is monochromatic.

Following the proof of \cite[{\sc Theorem}]{Loricolour} and using the
Preliminary Lemma, it is not hard to define a strategy for player~I
in such a way that if player~I follows this strategy, then for the
resulting partition $X$\;--\;which must belong to $\F$, since $\F$ is
a game filter\;--\;we get $s\seg X$ and $\nseg{s,X}{(|s|+k)}$ is
monochromatic. \qed

We have seen that every game filter has the segment-colouring
property. It seems that the reverse implication is unlikely, since a
strategy for player~I cannot be encoded by a real number, which makes
it hard (if not impossible) to prove that {\CH} implies the existence
of game filters. But on the other hand we know that scp-filters
always exist if we assume {\CH}.


\begin{thebibliography}{MmMm88a}

\bibitem[AdMo68]{AdMo68}
{\art{John\;W.~Addison and Yiannis\;N.~Moschovakis}
     {Some consequences of the axiom of definable determinateness}
     {Proceedings of the National Academy of Sciences U.S.A.}
     {59}
     {1968}
     {708--712}}

\bibitem[BaJu95]{BartoszynskiJudah}
{\book{Tomek~Bartoszy\'{n}ski and Haim~Judah}
    {Set Theory: On the Structure of the Real Line}
    {A.\,K.\,Peters}
    {Wellesley}
    {1995}}

\bibitem[BeKe84]{BeckerKechris}
{\samp{Howard\;S.~Becker and Alexander\; S.~Kechris}
     {Sets of ordinals constructible from trees and
      the third Victoria Delfino Problem}
     {Axiomatic Set Theory}
     {James\;E.~Baumgartner, Donald\;A.~Martin, Saharon~Shelah, {\sl Eds.}}
     {$[$Contemporary Mathematics 31$]$}
     {American Mathematical Society}
     {1984}
     {13--29}}

\bibitem[Br00${}_1$]{BrendlePrague}
{\samp{J\"org Brendle}
     {How small can the set of generics be?}
     {Logic Colloquium '98, Proceedings of the 1998 Association
      for Symbolic Logic European Summer Meeting, Prague,
      Czech Republic, 1998}
     {Samuel R.~Buss, Petr H\'ajek, Pavel Pudl\'ak, {\sl Eds.}}
     {$[$Lecture Notes in Logic 13$]$}
     {Springer}
     {2000}
     {92--109}}

\bibitem[Br00${}_2$]{BrendleMA}
{\art{J\"org~Brendle}
          {Martin's axiom and the dual distributivity number}
          {Mathematical Logic Quarterly}
          {46}
          {2000}
          {241--248}}

\bibitem[BrL\"o99]{BrLoe99}
{\art{J\"org~Brendle and Benedikt~L\"owe}
     {Solovay--type characterizations of forcing algebras}
     {Journal of Symbolic Logic}
     {64}
     {1999}
     {1307--1323}}

\bibitem[CaSi84]{CarlsonSimpson}
{\art{Timothy\;J.~Carlson and Steve\;G.~Simpson}
      {A dual form of Ramsey's Theorem}
      {Advances in Mathematics}
      {53}
      {1984}
      {265--290}}

\bibitem[El74]{Ellentuck}
{\art{Erik~Ellentuck}
     {A new proof that analytic sets are Ramsey}
     {Journal of Symbolic Logic}
     {39}
     {1974}
     {163--165}}

\bibitem[Halb98${}_1$]{Lorisym}
{\art{Lorenz~Halbeisen}
      {Symmetries between two Ramsey poperties}
      {Archive for Mathematical Logic}
      {37}
      {1998}
      {241--260}}

\bibitem[Halb98${}_2$]{Lorisha}
{\art{Lorenz~Halbeisen}
      {On shattering, splitting and reaping partitions}
      {Mathematical Logic Quarterly}
      {44}
      {1998}
      {123--134}}

\bibitem[Halb$\infty$]{Loricolour}
{\submitted{Lorenz~Halbeisen}
      {A Ramsey type theorem and its associated filters}}

\bibitem[HalbJu96]{LoriJudah}
{\art{Lorenz~Halbeisen and Haim~Judah}
     {Mathias absoluteness and the Ramsey property}
     {Journal of Symbolic Logic}
     {61}
     {1996}
     {177--193}}

\bibitem[HalbL\"o$\infty$]{luxor}
{\toappear{Lorenz~Halbeisen and Benedikt~L\"owe}
          {Ultrafilter spaces on the semilattice of partitions}
          {Topology and its Applications}}

\bibitem[HalJe63]{HalesJewett}
{\art{Alfred W.~Hales and Robert I.~Jewett}
     {Regularity and positional games}
     {Transactions of the American Mathematical Society}
     {106}
     {1963}
     {222--229}}

\bibitem[HarKe81]{HarringtonKechris}
{\art{Leo\;A.~Harrington and Alexander\;S.~Kechris}
     {On the determinacy of games on ordinals}
     {Annals of Mathematical Logic}
     {20}
     {1981}
     {109--154}}

\bibitem[Je78]{Jechbook}
{\book{Thomas~Jech}
    {Set Theory}
    {Academic Press}
    {San Diego}
    {1978}}

\bibitem[JuSh89]{JudahShelah.12}
{\art{Haim~Judah and Saharon~Shelah}
     {$\Delta^1_2$-sets of reals}
     {Annals of Pure and Applied Logic}
     {42}
     {1989}
     {207--223}}

\bibitem[JuSh93]{JudahShelah.13}
{\art{Haim~Judah and Saharon~Shelah}
     {$\Delta^1_3$-sets of reals}
     {Journal of Symbolic Logic}
     {58}
     {1993}
     {72--80}}

\bibitem[Kan94]{Kan94}
{\book{Akihiro~Kanamori}
      {The Higher Infinite}
      {$[$Perspectives in Mathe\-matical Logic$]$, Springer-Verlag}
      {Berlin Heidelberg}
      {1994}}

\bibitem[Kas83]{Kastanas}
{\art{Ilias\;G.~Kastanas}
    {On the Ramsey property for sets of reals}
    {Journal of Symbolic Logic}
    {48}
    {1983}
    {1035--1045}}

\bibitem[Ke95]{Kechris}
{\book{Alexander\;S.~Kechris}
    {Classical Descriptive Set Theory}
    {$[$Graduate Texts in Mathematics 156$]$, Springer-Verlag}
    {New York}
    {1995}}

\bibitem[KeMo72]{KeMo72}
{\art{Alexander\;S.~Kechris and Yiannis\;N.~Moschovakis}
    {Two theorems about projective sets}
    {Israel Journal of Mathematics}
    {12}
    {1972}
    {391--399}}

\bibitem[KeMo78]{KechrisMoschovakis}
{\samp{Alexander\;S.~Kechris and Yiannis\;N.~Moschovakis}
      {Notes on the Theory of Scales}
      {Cabal Seminar 76--77, Proceedings, Caltech--UCLA Logic Seminar 1976--77}
      {A.\;S.~Kechris and Y.\;N.~Moschovakis, {\sl Eds.}}
      {$[$Lecture Notes in Mathematics 689$]$, Springer-Verlag}
      {Berlin}
      {1978}
      {1--53}}

\bibitem[Ku83]{Kunen}
{\book{Kenneth~Kunen}
    {Set Theory, an Introduction to Independence Proofs}
    {$[$Studies in Logic and the Foundations of Mathematics 102$]$,
     North Holland}
    {Am\-st\-erdam}
    {1983}}
%
%


%

\bibitem[Mar68]{Ma68}
{\art{Donald\;A.~Martin}
     {The axiom of determinateness and reduction principles
      in the analytic hierarchy}
     {Bulletin of the American Mathematical Society}
     {74}
     {1968}
     {687--689}}

\bibitem[Mat77]{Mathias}
{\art{Adrian\;R.\;D.~Mathias}
     {Happy families}
     {Annals of Mathematical Logic}
     {12}
     {1977}
     {59--111}}

\bibitem[Mo71]{Mo71}
{\art{Yiannis\;N.~Moschovakis}
     {Uniformization in a playful universe}
     {Bulletin of the American Mathematical Society}
     {77}
     {1971}
     {731--736}}

\bibitem[Mo80]{Moschovakis}
{\book{Yiannis\;N.~Moschovakis}
      {Descriptive Set Theory}
      {$[$Studies in Logic and the Foundations of
       Mathematics 100$]$, North-Holland}
      {Amsterdam}
      {1980}}

\bibitem[Pl86]{Plewik}
{\art{Szymon~Plewik}
     {On completely Ramsey sets}
     {Fundamenta Mathematicae}
     {127}
     {1986}
     {127--132}}

\bibitem[Rai84]{Raisonnier}
{\art{Jean~Raisonnier}
     {A mathematical proof of S.~Shelah's Theorem on the
      Measure Problem and related results}
     {Israel Journal of Mathematics}
     {48}
     {1984}
     {48--56}}

\bibitem[Ram29]{Ramsey}
{\art{Frank\;P.~Ramsey}
     {On a problem of formal logic}
     {Proceedings of the London Mathematical Society, Ser.\,II}
     {30}
     {1929}
     {264--286}}

\bibitem[St95]{PWOIM}
{\art{John\;R.~Steel}
     {Projectively well--ordered inner models}
     {Annals of Pure and Applied Logic}
     {74}
     {1995}
     {77--104}}

\end{thebibliography}
\end{document}